\documentclass[preprint]{imsart}

\usepackage[T1]{fontenc}

\RequirePackage{amsthm,amsmath,amssymb,amsfonts,calrsfs}
\RequirePackage[colorlinks,citecolor=blue,urlcolor=blue]{hyperref}
\RequirePackage{hypernat}

\usepackage{epsfig}
\usepackage{pst-plot}
\usepackage{multirow}

\startlocaldefs

\numberwithin{equation}{section}

\newtheorem{proposition}{Proposition}[section]
\newtheorem{theorem}[proposition]{Theorem}

\newtheorem{example}[proposition]{Example}

\newcommand{\R}{{\mathbb R}}

\newcommand{\N}{{\mathbb N}}
\newcommand{\Z}{{\mathbb Z}}

\newcommand{\Pareto}{\operatorname{Pareto}}

\newcommand{\harm}{\operatorname{harm}}

\renewcommand{\le}{\leqslant}
\renewcommand{\ge}{\geqslant}

\newcommand{\norm}[1]{\|#1\|}

\newcommand{\E}{\operatorname{E}}

\newcommand{\biindice}[3]%
{

\begin{array}[t]{c}
#1\\
{\scriptstyle #2}\\
{\scriptstyle #3}
\end{array}

}

\endlocaldefs

\begin{document}

\begin{frontmatter}
\title{Maxima of moving maxima of continuous functions\protect\thanksref{}}
\runtitle{CM3 processes}

\author{\fnms{Thomas} \snm{Meinguet}\thanksref{t1}\ead[label=e1]{thomas.meinguet@uclouvain.be}}

\thankstext{t1}{Research supported by IAP research network grant nr.\ P6/03 of the Belgian government (Belgian Science Policy) and by contract nr.\ 07/12/002 of the Projet d'Actions de Recherche Concert\'ees of the Communaut\'e fran\c{c}aise de Belgique, granted by the Acad\'emie universitaire Louvain.}

\runauthor{T. Meinguet}

\affiliation{Universit\'{e} catholique de Louvain}

\address{Universit\'e catholique de Louvain, Institut de statistique\\
Voie du Roman Pays 20, B-1348 Louvain-la-Neuve, Belgium\\
\printead{e1}}

\begin{abstract}
Maxima of moving maxima of continuous functions (CM3) are max-stable processes aimed at modeling extremes of continuous phenomena over time. They are defined as Smith and Weissman's M4 processes with continuous functions rather than vectors. After standardization of the margins of the observed process into unit-Fréchet, CM3 processes can model the remaining spatio-temporal dependence structure.

CM3 processes have the property of joint regular variation. The spectral processes from this class admit particularly simple expressions. Furthermore, depending on the speed with which the parameter functions tend toward zero, CM3 processes fulfill the finite-cluster condition and the strong mixing condition. For instance, these three properties put together have implications for the expression of the extremal index. 

A method for fitting a CM3 to data is investigated. The first step is to estimate the length of the temporal dependence. Then, by selecting a suitable number of blocks of extremes of this length, clustering algorithms are used to estimate the total number of different profiles. The number of parameter functions to retrieve is equal to the product of these two numbers. They are estimated thanks to the output of the partitioning algorithms in the previous step. The full procedure only requires one parameter which is the range of variation allowed among the different profiles. The dissimilarity between the original CM3 and the estimated version is evaluated by means of the Hausdorff distance between the graphs of the parameter functions.
\end{abstract}



\end{frontmatter}

\section{Introduction}

\emph{Maxima of moving maxima of continuous functions} (CM3) are the analogue of Smith and Weissman's M4 processes \cite{SW96} with continuous functions rather than vectors.
Let $a_i^{(j)}$ ($i\in\Z_+$, $j\in\Z$) be strictly positive, real, continuous functions on a compact domain of $\R^q$, say $[0,1]^q$. The functions $a_i^{(j)}$ are the \emph{parameter functions}. They are assumed to satisfy, for every $x\in[0,1]^q$, the equality $\sum_{j\in\Z}\sum_{i\geqslant0} a_i^{(j)}(x)=1.$
A CM3 process $(X_t)_{t\in\Z}$ is defined by the expression
$$ X_t(x) = \sup_{j\in\Z} \sup_{i \geqslant 0} a_i^{(j)}(x) Z_{t-i}^{(j)}\qquad (x\in[0,1]^q),$$
where the innovations $Z_i^{(j)}$ ($i\in\Z$, $j\in\Z$) are independent and identically distributed unit-Fr\'echet random variables,
i.e.~$P(Z_t\leqslant z)=\exp(-1/z)$ for $z>0$.

The fact that, given real numbers $\xi_i>0$ ($i\in\N$) such that $\xi_1 + \xi_2 + \ldots =1$, the distribution of
$ \max(\xi_1 Z_1,\xi_2 Z_2,\ldots)$ stays unit-Fr\'echet
implies that $X_t$ has unit-Fr\'echet margins. However the transformation from $(Z_t)_{t\in\Z}$ to $(X_t)_{t\in\Z}$ induces a dependence
structure in time and space. Extremes appear in temporal clusters and, at time $t$, a large value for $X_t$ at location $x$ causes
large values at other locations. From this fact, CM3 processes are able to model a wide range of spatio-temporal dependences. The first part of this paper is a study of some properties: spectral process, strong mixing condition, finite-cluster condition and extremal index.

The second objective of this paper is to fit CM3 processes to samples with measurement errors.
For that purpose, CM3 will be discretized into M4 of dimension $D$ selecting $D$ points $x_d$ $(1 \leqslant d \leqslant D)$ in the domain.
It will be also assumed that $0\leqslant i< K$ and $1\leqslant j \leqslant L$ for finite constants $K$ and $L$.
The practical model studied is thus
$$ X_t(d) = \max_{1 \leqslant j \leqslant L} \max_{0 \leqslant i < K} a_i^{(j)}(x_d) Z_{t-i}^{(j)}+\varepsilon_t(x_d)\qquad (1 \leqslant d \leqslant D)$$
where $\varepsilon_t(x_d)$ are independent $N(0,\sigma^2)$ random variables.
The parameter $K$ is the length of the temporal dependence and $L$ is the total number of reproducible patterns that we can observe up to a multiplicative constant in the process.

\begin{figure}[here!]
\begin{center}
\includegraphics[scale=0.325]{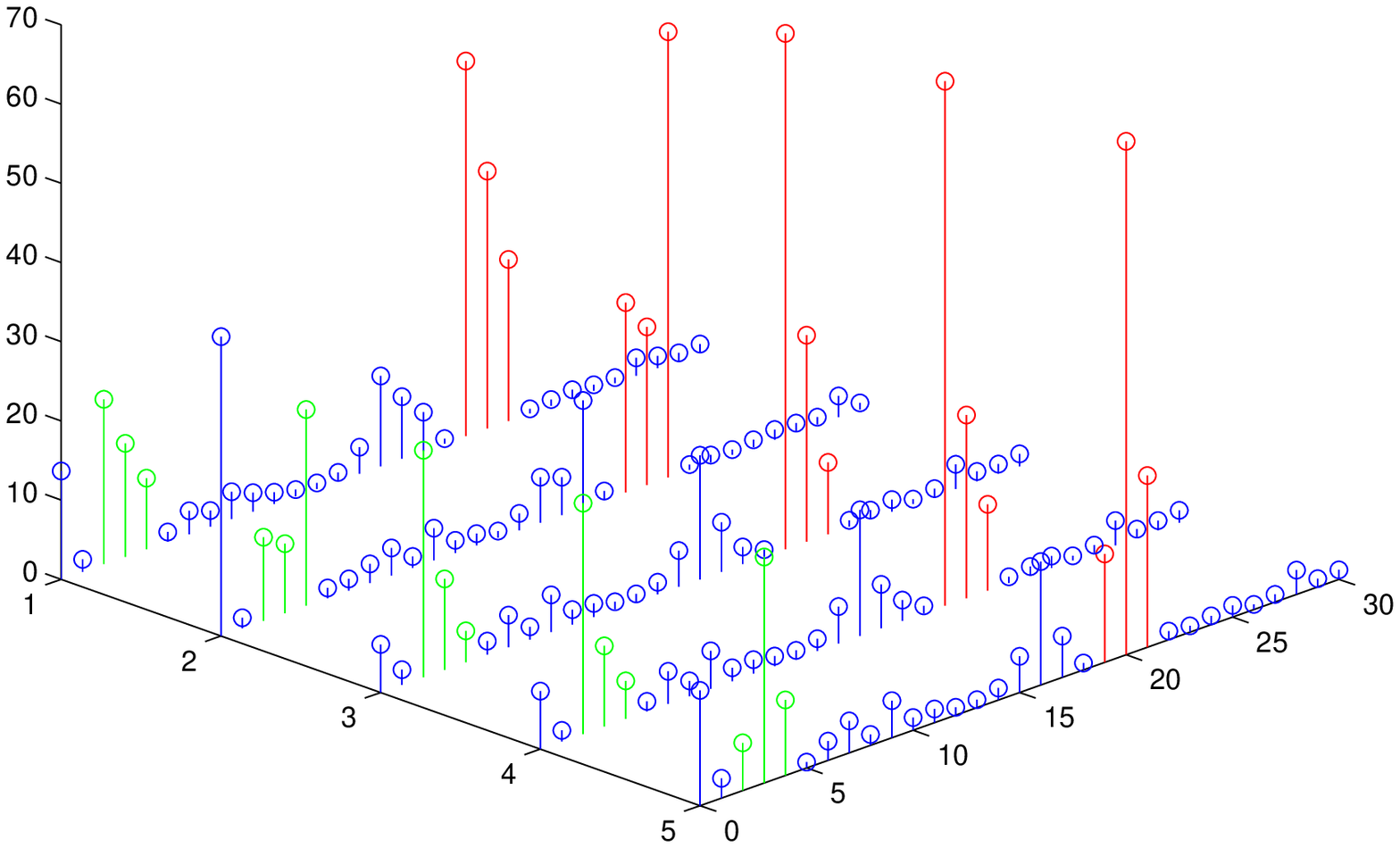}
\includegraphics[scale=0.325]{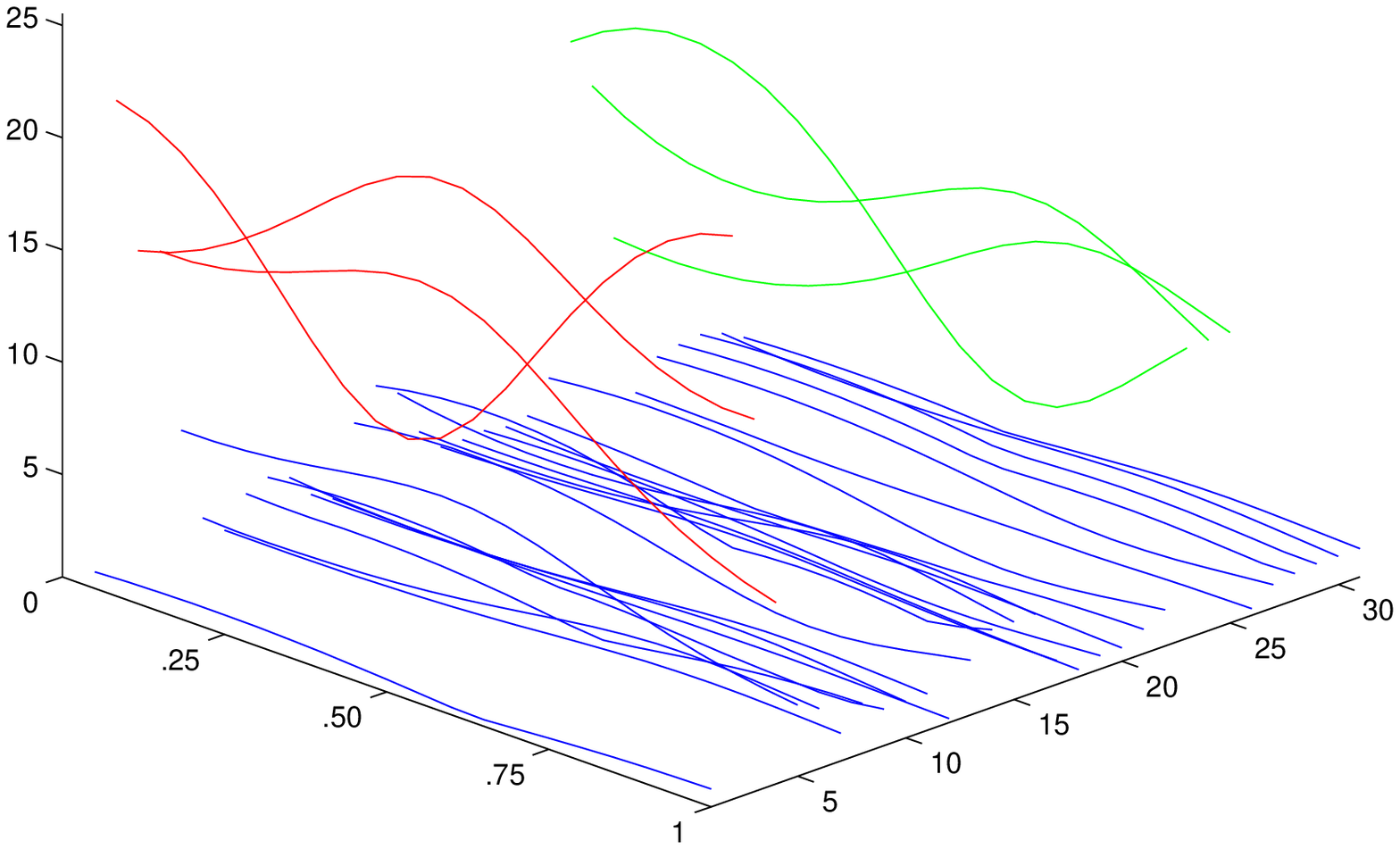}
\begin{pspicture}(-7,0)(7,0)
\psset{xunit=1cm,yunit=1cm,runit=1cm}
\uput[r](-2.2,.7){\footnotesize{$t$}} \uput[r](3.7,.7){\footnotesize{$t$}}
\uput[r](-6,.7){\footnotesize{$d$}} \uput[r](-.3,.7){\footnotesize{$x$}}
\uput[r](-6.3,3.7){\footnotesize{$X_t(d)$}} \uput[r](-.6,3.7){\footnotesize{$X_t(x)$}}
\end{pspicture}
\caption{M4 with $D=5$ on the left, CM3 on $[0,1]$ on the right ($K=3$, $L=2$).}
\label{M-CM}
\end{center}
\end{figure}

Figure~\ref{M-CM} shows a realization of a CM3 plotted versus a M4.

In Section~\ref{SectionCM3}, a coherent set of properties for CM3 is established.
The motivation is similar as in \cite{SEG06,SMI03} but now for random continuous functions.
Theorem~\ref{CM3process} is the joint regular variation of those processes,
a concept extended to Banach spaces in \cite{MS10}.
The spectral process of a CM3 has a discrete distribution, given by the theorem.
Next, depending on the speed with which the parameter functions $a_i^{(j)}$ tend
toward zero, Theorem~\ref{conditionCThm} yields
the finite-cluster condition and Theorem~\ref{conditionSMThm} yields the strong
mixing condition. These three properties together also have specific implications,
for instance the inverse of the extremal index $\theta$ becomes the expected size of
clusters of extremes in the sense of \cite{SEG05}.

CM3 processes are also examples of max-stable random fields \cite{DM08,DHF06}: for every finite space-time subset $A\times T\subset[0,1]^q\times\Z$, the random vector $(X_t(x))_{x\in A,t\in T}$ has a multivariate extreme value distribution. This property of M4 is inherent to CM3 since the law of a continuous random field is characterized by its finite dimensional distributions that are M4
according to Example~\ref{ExampleM4}.

Section~\ref{blocksprofile} is a preparation for the estimation of the parameter functions.
Extremes will play a central role in identifying the recursive patterns and their relative frequencies.
So we need to study the probabilistic properties of the blocks of extremes that can be observed in CM3.
The harmonic mean makes convenient the expressions of the frequencies of the reproducible patterns
that can be observed.
 
In Section~\ref{Estimation}, we suggest and compare empirical methods to
estimate $K$, $L$ and the parameter functions without assumptions on
these objects. It is a complement to \cite{Z04}, where the case $L=1$
has only been treated, and to \cite{Z08}, where assumptions are made on
the parameters.

This study is designed to improve the statistical analyses of
extreme events, as done in \cite{SUV09,SD10} for instance.

\section{Definition and properties} \label{SectionCM3}

Choose a nonempty compact domain of $\R^q$. To not multiply the notations this compact will be taken to be $[0,1]^q$.
Given an array $Z_i^{(j)}$ ($i\in\Z_+$, $j\in\Z$) of independent unit-Fr\'{e}chet random variables, if $a_i^{(j)}:[0,1]^q\to\R_+^*$
are deterministic strictly positive continuous functions, a \emph{CM3 process} is a stochastic process defined by
\begin{equation} \label{CM3def}
X_t(x)= \sup_{j\in\Z} \sup_{i\geqslant0} (a_i^{(j)}(x) Z_{t-i}^{(j)}).
\end{equation}
If furthermore $$ (\forall x\in[0,1]^q):\sum_{j\in\Z}\sum_{i\geqslant0}
a_i^{(j)}(x)=1$$ we say that $(X_t)_{t\in\Z}$ is a \emph{standard CM3 process}.

The first result is an imperative condition before any use of CM3 processes. Recall that the sup-norm of a function $f:[0,1]^q\to\R$
is $ \| f \|_{\infty} = \sup_{x\in[0,1]^q}|f(x)|$ and that this supremum is achieved.

\begin{proposition} \label{CM3processLem}
If
\begin{equation} \label{finiteSum}
\sum_{j \in\Z}\sum_{i\geqslant0} \| a_i^{(j)} \|_\infty < \infty,
\end{equation}
then, for every $t\in\Z$, $X_t$ in \eqref{CM3def} is a random element in $\mathcal{C}([0,1]^q,\R_+)$
and the process $(X_t)_{t\in\Z}$ is stationary.
\end{proposition}
The proofs of the results of this section are relegated to Appendix~\ref{proofs}.

\begin{example} \label{ExampleM4}
If $(X_t)_{t\in\Z}$ is a standard CM3 with $0\leqslant i < K$ and $1 \leqslant j \leqslant L$,
then for each $x_1,\ldots,x_D\in\R^q$, the process $(X_t(x_1),\ldots,X_t(x_D))_{t\in\Z}$ is a standard M4. Under \eqref{finiteSum}
$(\|X_t\|_\infty)_{t\in\Z}$ is a non-standard M3. In both cases $0\leqslant i < K$ and $1\leqslant j \leqslant L$.
This is helpful for the estimation of $K$ and $L$ if $x_1,\ldots,x_D$ are far enough, see \ref{ChooseExt}.
\end{example}

A CM3 process is an example of a jointly regularly varying time series. In particular, there exists
a process $\left(\Theta_t\right)_{t\in\mathbb{Z}}$ in $\mathcal{C}([0,1]^q,\R_+)$, called \emph{spectral process} which is
the limit in distribution, as $x\to\infty$,
$$ \mathcal{L}\left(\left.\left(X_t/\|X_0\|_\infty\right)_{t\in\mathbb{Z}}\ \right|\ \| X_0 \|\geqslant x \right)\\
 \xrightarrow{d} \mathcal{L}\left(\left(\Theta_t\right)_{t\in\mathbb{Z}}\right) $$
in the proper product space. According to \cite{MS10}, this process captures all aspects of extremal dependence,
both within space and over time.

\begin{theorem} \label{CM3process}
Setting $a_i^{(j)}=0$ if $i<0$, under condition \eqref{finiteSum}, a CM3 process $(X_t)_{t\in Z}$ is jointly regularly varying with index $\alpha=1$ and spectral process
$$
(\Theta_{-s},\ldots,\Theta_{t}) \overset{d}{=} \left(\frac{a_{-s+I}^{(J)}}{\| a_{I}^{(J)}\|_\infty} ,\ldots, \frac{a_{t+I}^{(J)}}{\| a_{I}^{(J)}\|_\infty}\right),
$$
where $(I,J)$ is a random vector on $\Z_+\times\Z$ having distribution
$$ P[(I,J)=(i,j)]= \frac{\|a_i^{(j)}\|_\infty}{\sum_{l\in\Z}\sum_{k\geqslant 0} \|a_k^{(l)}\|_\infty},\ i\in\Z_+,\ j\in\Z.  $$
\end{theorem}

All CM3 processes satisfying \eqref{finiteSum} also satisfy the finite-cluster condition. This property prevents a sequence of extremes occurring in a CM3 from being infinite over time even if $K=+\infty$ or $L=+\infty$.

\begin{theorem} \label{conditionCThm}
Under condition \eqref{finiteSum}, a CM3 process $(X_t)_{t\in Z}$ satisfies the finite-cluster condition:
there exists $(r_n)_{n\in\N}$ with $r_n\to\infty$ and $r_n/n\to 0$ such that
\begin{equation} \label{conditionC} \tag{C}
\lim_{m\to\infty} \limsup_{n\to\infty} P(\max_{m\leqslant|t|\leqslant r_n}\|X_t\|_\infty>n\ |\ \|X_0\|_\infty>n )=0.
\end{equation}
\end{theorem}

Together with the finite-cluster condition, the strong mixing property leads to nice properties. To obtain the strong mixing property
a sufficient condition is
\begin{equation} \label{fastConv}
\sum_{j \in\Z}\sum_{i\geqslant0} i \| a_i^{(j)} \|_\infty < \infty.
\end{equation}
Note that \eqref{finiteSum} and \eqref{fastConv} are trivial whenever $K<+\infty$ and $L<+\infty$.

\begin{theorem} \label{conditionSMThm}
Under condition \eqref{fastConv}, a CM3 process $(X_t)_{t\in Z}$
satisfies the strong mixing condition:
\begin{equation} \label{strongMixing} \tag{M}
\lim_{m\to\infty} \biindice{\sup}{A\in\sigma(-\infty,-m)}{B\in\sigma(m,\infty)} |P(A\cap B)-P(A)P(B)|=0
\end{equation}
where $\sigma(r,s)$ is the $\sigma$-field generated by $\{X_t\, |\, r\leqslant t \leqslant s\}$.
\end{theorem}

If $(X_t)_{t\in\Z}$ is a regularly varying time series with index 1, the extremal index $\theta$ of the univariate time series $(\|X_t\|_\infty)_{t\in\Z}$ is defined as the quantity between $0$ and $1$ such that
$$ P(\max_{1\leqslant t \leqslant n}\|X_t\|_\infty\leqslant nx) \to e^{-\theta /x} $$
as $n\to\infty$.
The extremal index of a CM3 process is the following.
\begin{proposition} \label{extIndex}
Under condition \eqref{finiteSum}, if $(X_t)_{t\in Z}$ is a CM3 process, the extremal index of $(\|X_t\|_\infty)_{t\in\Z}$ is
$$ \theta=\frac{1}{\displaystyle \sum_{j\in\Z}\sum_{i\geqslant0}\|a_i^{(j)}\|_\infty}\ \sum_{j\in\Z}\max_{i\geqslant0}\|a_i^{(j)}\|_\infty. $$
\end{proposition}

Once conditions \eqref{conditionC} and \eqref{strongMixing} are satisfied, which is the case under \eqref{fastConv} by Theorem~\ref{conditionCThm} and Theorem~\ref{conditionSMThm}, there are further characterizations of the extremal index such that
\begin{equation} \label{exiexp}
\begin{array}{rcl}
\theta &=& \displaystyle  \lim_{t\to\infty} \lim_{x\to\infty} P(\max_{i=1,\ldots,t} \|X_i\|_\infty \leqslant x\ |\ \|X_0\|_\infty>x )\\
\end{array}
\end{equation}
and, in this case, $1/\theta$ is the expected size of clusters of extremes in the sense of \cite{SEG05}, which is recalled as follows.
Let $u_n\to\infty$ be a
thresholding
sequence
and $r_n\to\infty$ be such that the expected number of exceedances in a sample of size $r_n$ tends toward 0:
$$ E\left[ \sum_{i=1}^{r_n} 1\{\|X_i\|_\infty>u_n\} \right]=r_nP(\|X_1\|_\infty>u_n)\to 0.$$
Then, denoting $M_n:=\max\{\|X_1\|_\infty,\ldots,\|X_{n}\|_\infty\}$, under \eqref{conditionC} and \eqref{strongMixing}, we have
$$ E\left[ \left.\sum_{i=1}^{r_n} 1\{\|X_i\|_\infty>u_n\}\ \right|\ M_{r_n}>u_n \right]
= \frac{r_nP(\|X_1\|_\infty>u_n)}{P( M_{r_n}>u_n)}\to\frac{1}{\theta}$$
as $n\to\infty$.

\section{Block profiles} \label{blocksprofile}

In this section we study further probabilistic features of CM3 processes in order to build a method to estimate the parameter functions
$a_i^{(j)}$ in the case $0\leqslant i< K$ and $1\leqslant j \leqslant L$.
The theoretical model for the rest of the paper is thus
\begin{equation} \label{theoreticalModel}
X_t(x) = \max_{1 \leqslant j \leqslant L} \max_{0 \leqslant i < K} a_i^{(j)}(x) Z_{t-i}^{(j)}\qquad (x\in[0,1]^q)
\end{equation}

The estimation method suggested in this paper is based on the fact that a large value of $Z_{i}^{(j)}$ causes large values of $X_t$ for
$i\leqslant t < i+K$ and the possibility to have in this case
\begin{equation} \label{mainEvent}(X_t,\ldots,X_{t+K-1})=Z_{i}^{(j)}(a_0^{(j)},\ldots,a_{K-1}^{(j)}).\end{equation}
By ``block profile'' we mean a sequence $(X_t,\ldots,X_{t+K-1})$ satisfying  \eqref{mainEvent} for some\linebreak$1\leqslant j \leqslant L$. The corresponding sequence of functions $(a_0^{(j)},\ldots,a_{K-1}^{(j)})$ will be called ``profile'' or ``pattern''.

In \ref{relfreq} we compute the probability of the events \eqref{mainEvent} and their frequencies of occurrence for the different values of $j$. In \ref{correlation} we have a brief look at the correlation between all the possible blocks of length $K$ available in a sample.
They are not independent if they overlap. In \ref{samplesize} we give the needed sample size to expect that an event of the form $\eqref{mainEvent}$ realizes at least once.

To compute the exact values, the knowledge of the parameter functions is needed, which is particular not the case in the estimation. This is the reason why we also give lower and upper bounds for the true values. These bounds only depend on a unique parameter $C$, which is the maximal variation among the parameter functions $a_i^{(j)}$.

\subsection{Relative frequencies} \label{relfreq}

To recover the functions $a_i^{(j)}$ from a sample of size $T$,
the first step is to understand how \eqref{theoreticalModel} works.
Consider as an example a simple situation when $K=3$ and $L=2$.
A finite number of functions $a_i^{(j)}$ is uniformly bounded below by a positive constant
since they are strictly positive. Thus, if for instance $Z_{1}^{(2)}$ is large enough,
the value of $(X_1,\ldots,X_K)$ at a given position $x\in[0,1]^q$ is
\begin{equation} \label{matrices}
\begin{array}{c}
X_1(x)=\max\left(
\begin{array}{lll}
a_2^{(1)}(x)Z_{-1}^{(1)}\ \ & a_1^{(1)}(x)Z_{0}^{(1)}\ \ & a_0^{(1)}(x)Z_{1}^{(1)}\\
a_2^{(2)}(x)Z_{-1}^{(2)} & a_1^{(2)}(x)Z_{0}^{(2)} & a_0^{(2)}(x)\mathbf{Z_{1}^{(2)}}\\
 \end{array}\right)= a_0^{(2)}(x)Z_{1}^{(2)}

 \\

\rule{0cm}{1cm}

X_2(x)=\max\left(
\begin{array}{lll}
a_2^{(1)}(x)Z_{0}^{(1)}\ \ & a_1^{(1)}(x)Z_{1}^{(1)}\ \ & a_0^{(1)}(x)Z_{2}^{(1)}\\
a_2^{(2)}(x)Z_{0}^{(2)} & a_1^{(2)}(x)\mathbf{Z_{1}^{(2)}} & a_0^{(2)}(x)Z_{2}^{(2)}\\
 \end{array}\right)= a_1^{(2)}(x)Z_{1}^{(2)}

 \\

\rule{0cm}{1cm}

X_3(x)=\max\left(
\begin{array}{lll}
a_2^{(1)}(x)Z_{1}^{(1)}\ \ & a_1^{(1)}(x)Z_{2}^{(1)}\ \ & a_0^{(1)}(x)Z_{3}^{(1)}\\
a_2^{(2)}(x)\mathbf{Z_{1}^{(2)}} & a_1^{(2)}(x)Z_{2}^{(2)} & a_0^{(2)}(x)Z_{3}^{(2)}\\
 \end{array}\right)= a_2^{(2)}(x)Z_{1}^{(2)}

\end{array}
\end{equation}
so that the second pattern appears: $$(X_1,X_2,X_3)(x)=(a_0^{(2)},a_1^{(2)},a_2^{(2)})(x)Z_{1}^{(2)}.$$
How likely is this kind of events to occur? To compute their probabilities, first remark that \eqref{matrices}
is equivalent to
\begin{equation} \label{intersection}
\left\{\begin{array}{ll}
\displaystyle Z_{-1}^{(1)} \leqslant \min_x[\frac{a_0^{(2)}(x)}{a_2^{(1)}(x)} ]Z_{1}^{(2)}
&  \displaystyle Z_{-1}^{(2)} \leqslant \min_x[\frac{a_0^{(2)}(x)}{a_2^{(2)}(x)}] Z_{1}^{(2)} \\
\displaystyle Z_{0}^{(1)} \leqslant \min_x[\frac{a_0^{(2)}(x)}{a_1^{(1)}(x)},\frac{a_1^{(2)}(x)}{a_2^{(1)}(x)} ] Z_{1}^{(2)}
&  \displaystyle Z_{0}^{(2)} \leqslant \min_x[\frac{a_0^{(2)}(x)}{a_1^{(2)}(x)},\frac{a_1^{(2)}(x)}{a_2^{(2)}(x)} ] Z_{1}^{(2)} \\
\displaystyle Z_{1}^{(1)} \leqslant \min_x[\frac{a_0^{(2)}(x)}{a_0^{(1)}(x)},\frac{a_1^{(2)}(x)}{a_1^{(1)}(x)},\frac{a_2^{(2)}(x)}{a_2^{(1)}(x)} ] Z_{1}^{(2)}
& \qquad\qquad\qquad -\\
\displaystyle Z_{2}^{(1)} \leqslant \min_x[\frac{a_1^{(2)}(x)}{a_0^{(1)}(x)},\frac{a_2^{(2)}(x)}{a_1^{(1)}(x)} ] Z_{1}^{(2)}
&  \displaystyle Z_{2}^{(2)} \leqslant \min_x[\frac{a_1^{(2)}(x)}{a_0^{(2)}(x)},\frac{a_2^{(2)}(x)}{a_1^{(2)}(x)} ] Z_{1}^{(2)} \\
\displaystyle Z_{3}^{(1)} \leqslant \min_x[\frac{a_2^{(2)}(x)}{a_0^{(1)}(x)} ]Z_{1}^{(2)}
&  \displaystyle Z_{3}^{(2)} \leqslant \min_x[\frac{a_2^{(2)}(x)}{a_0^{(2)}(x)} ]Z_{1}^{(2)} \\
\end{array}\right.
\end{equation}
For the general case, let $A(l^\ast)$ be the event \eqref{mainEvent} with fixed $j=l^\ast$:
$$A(l^\ast)=\{\forall x\in[0,1]^q:(X_t,\ldots,X_{t+K-1})(x)=Z_t^{(l^\ast)}(a_0^{(l^\ast)},\ldots,a_{K-1}^{(l^\ast)})(x)\},$$
i.e.~$A(l^\ast)$ is the event for a $K$-block starting a time $t$ to be a block profile of type~${l^\ast}$.
Generalizing \eqref{intersection} shows that the event $A(l^\ast)$ is the intersection of $(2K-1)L-1$ conditions
involving the random variables $Z_i^{(j)}$ for $t-K+1\leqslant i \leqslant t+K-1$.
Remembering that the density of $Z$ is $f_Z(z)=z^{-2}\exp(-z^{-1})$, the probability $p^{(l^\ast)}$ of $A(l^\ast)$ is
\begin{equation}\label{bigintegral}
\begin{array}{rcl}
\displaystyle p^{(l^\ast)} &=& \displaystyle \int_0^\infty  \prod_{l=1}^L P(Z\leqslant \min_{x}[\frac{a_0^{(l^\ast)}(x)}{a_{K-1}^{(l)}(x)}] z )\\
&& \displaystyle \prod_{l=1}^L P(Z\leqslant \min_{x}[\frac{a_0^{(l^\ast)}(x)}{a_{K-2}^{(l)}(x)},\frac{a_1^{(l^\ast)}(x)}{a_{K-1}^{(l)}(x)}] z )\\
&& \ldots\\
&& \displaystyle \!\!\!\!\!\biindice{\displaystyle \prod^L}{l=1}{\ l\neq l^\ast} P(Z\leqslant \min_{x}[\frac{a_0^{(l^\ast)}(x)}{a_{0}^{(l)}(x)},\ldots,\frac{a_{K-1}^{(l^\ast)}(x)}{a_{K-1}^{(l)}(x)}] z )\\
&& \ldots\\
&& \displaystyle \prod_{l=1}^L P(Z\leqslant \min_{x}[\frac{a_{K-2}^{(l^\ast)}(x)}{a_0^{(l)}(x)},\frac{a_{K-1}^{(l^\ast)}(x)}{a_1^{(l)}(x)}] z )\\
&& \displaystyle \prod_{l=1}^L P(Z\leqslant \min_{x}[\frac{a_{K-1}^{(l^\ast)}(x)}{a_0^{(l)}(x)}] z ) z^{-2}\exp(-z^{-1})dz\\
\end{array}
\end{equation}
Denoting $m_k^{(l;l^\ast)}$ the minimum written on line $k$ of \eqref{bigintegral} and by extension $m_k^{(l^\ast;l^\ast)}=1$, the value of $p^{(l^\ast)}$ is
\begin{equation} \label{profprob}
p^{(l^\ast)} = \frac{1}{\displaystyle 1+\sum_{k=1}^{2K-1}  \biindice{\displaystyle\sum^L}{l=1}{l\neq l^\ast\textrm{ if }k=K}\frac{1}{m_k^{(l;l^\ast)}}}=\frac{\harm(m_1^{(1;l^\ast)},\ldots,m_{2K-1}^{(L;l^\ast)})}{(2K-1)L}\end{equation}
where $\harm$ is the harmonic mean of the $(2K-1)L$ minima.

It is thus possible to compute $p^{(l^\ast)}$ exactly given the parameter functions. If the parameter functions $a_i^{(j)}$ satisfy
\begin{equation}\label{bounds}
\frac{1}{C^{(l^\ast)}} \leqslant \frac{a_k^{(l^\ast)}(x)}{a_{k'}^{(l)}(x)} \leqslant C^{(l^\ast)}\end{equation}
for all $x,k,k',l,l^\ast$, then the probability $p^{(l^\ast)}$ of success to reveal $(a_0^{(l^\ast)},\ldots,a_{K-1}^{(l^\ast)})$ by picking up a random block satisfies
\begin{equation}
\underline{p}^{(l^\ast)}:=\frac{1}{C^{(l^\ast)}(2K-1)L}\leqslant p^{(l^\ast)}\leqslant\frac{C^{(l^\ast)}}{(2K-1)L}=:\overline{p}^{(l^\ast)} . \end{equation}
The harmonic mean being more sensitive to small values, the lower bound is actually closer to the exact probability.

Under the knowledge of $K$, the probability that a random $K$-block is a profile differs from pattern to pattern. But under the control condition \eqref{bounds}, if all $C^{(l^\ast)}$ are themselves bounded above by a common constant $C$,
the probability $p=p^{(1)}+\ldots+p^{(L)}$ for a random block to be any profile can be estimated by
\begin{equation}
\underline{p}:=\frac{1}{C(2K-1)} \leqslant p.
\end{equation}
which has the remarkable property not to depend neither on $L$ nor on the dimension of the ambient space.

As an illustration, Table \ref{TableBlocs} shows the number of found block profiles found versus their expectations,
knowing and without knowing the parameter functions for five simulations of \eqref{practicalModel}.
The different patterns are split in columns.

\begin{table}[here!]
\begin{center}
\begin{tabular}{|c|c|c|c|c|c||c||c|c|c|c|c|c|}
\hline
\multicolumn{13}{|c|}{Simulation with parameters $C=5$, $D=20$, $K=5$, $L=5$ and $T=5000$}\\
\hline \hline
\multicolumn{6}{|c||}{Expected value ($p^{(l^\ast)}T$)} & $\underline{p}T$ & \multicolumn{6}{|c|}{Really found}\\
\hline
\#1 & \#2 & \#3 & \#4 & \#5 & total & total & total & \#1 & \#2 & \#3 & \#4 & \#5 \\
\hline
31 & 31 & 30  & 34  & 33 & 159 & 135 & 147 & 28 & 20 & 25  & 30  & 44\\
31 & 32  & 33 & 29 & 36 & 161 & 135 & 159 & 30 & 20 & 43 & 36 & 30\\
35 & 35 & 31 & 31 & 33 & 166 & 135 & 166 & 31 & 35 & 24 & 30 & 46\\
32 & 32 & 36 & 33 & 33 & 165 & 135 & 160 & 33 & 26 & 43 & 27 & 31\\
30 & 30 & 30 & 31 & 33 & 154 & 135 & 143 & 31 & 28 & 32 & 19 & 33\\
\hline
\end{tabular}
\end{center}
\caption{Estimation of the number of block profiles versus values obtained in simulation.\label{TableBlocs}}
\end{table}

\subsection{Correlation} \label{correlation}

As we have seen in paragraph \ref{relfreq}, we need $(2K-1)L/\harm(m_1^{(1;l^\ast)},\ldots,m_{2K-1}^{(L;l^\ast)})$
independent random blocs from the series $(X_t)_{t\in\Z}$ to expect that at least one is proportioned like the ${l^\ast}^\textrm{th}$ profile.
Practically, given a chain $X_1,\ldots,X_T$ with $T$ observations, we have $T-K+1$ dependent blocks of length $K$.
The main pieces of information about the dependence structure between these blocks can be summarized in the following way.
\begin{enumerate}
\item[i)] If a $K$-block $(X_t,\ldots,X_{t+K-1})$ is a block profile, it does not overlap with another block profile.
\item[ii)] Given that the $K$-block $(X_t,\ldots,X_{t+K-1})$ is not a profile, the probability that one of the $K-1$ next $K$-blocks $(X_{t+1},\ldots,X_{t+K}),\ldots$ is a profile is higher.
\item[iii)] If consecutive blocks are not profiles, the probability that the next one is a profile stops increasing after $K$ non-profile blocks.
\end{enumerate}

To see i), for instance, have a look at the third matrix in \eqref{matrices}. If $(X_1,\ldots,X_3)$ is like the second profile, then in particular $$a_2^{(2)}Z_{1}^{(2)}>a_0^{(1)}Z_{3}^{(1)}.$$ But for $(X_3,\ldots,X_5)$ to be like the first profile, we must have $$a_2^{(2)}Z_{1}^{(2)}<a_0^{(1)}Z_{3}^{(1)}.$$ Then $(X_3,\ldots,X_5)$ cannot be proportioned like the first profile. Proceed similarly for any two non-disjoint blocks and any two different profiles.

To see ii), if for instance $(X_1,\ldots,X_3)$ in \eqref{intersection} is not a profile, that means that at least one of the 15 inequalities is not satisfied, although we do not know precisely how many. Some of the reverse inequalities lie in the conditions for the $K-1$ next blocks to be profile and some do not. To get the exact incidence, we need to condition on the number of inequalities not satisfied in \eqref{intersection} and whether or not they participate in the conditions for the next blocks to be profile. In any case: the probability that the next blocks are profile increases.

To see iii), simply remark that the process is $K$-dependent.

\subsection{Sample size} \label{samplesize}

As a consequence of paragraph \ref{correlation}, the expected number of profiles of type $l$ is greater than unity in a chain of length at least
\begin{equation}T=K-1+\frac{(2K-1)L}{\harm(m_1^{(1;l)},\ldots,m_{2K-1}^{(L;l)})}.\end{equation}
Given an upper bound $C$ on all the $C^{(l^\ast)}$ in \eqref{bounds}, the minimum sample size needed to expect at least $M$ repetitions of a particular profile in the chain is
\begin{equation}T=K-1+CM(2K-1)\end{equation}
if we do not know the parameter functions but only $C$.

\section{Estimation}  \label{Estimation}

The estimation methodology for the parameter functions $a_i^{(j)}$ starts from a discretization at $D$ points $x_d$ $(1 \leqslant d \leqslant D)$ of the domain. CM3 processes from \eqref{theoreticalModel} are seen as a high-dimensional M4.
Furthermore we may want to consider independent and normally distributed errors with variance $\sigma^2$ at each measurement point.
Thus the ``practical model'' studied in this section is
\begin{equation} \label{practicalModel}
X_t(d) = \max_{1 \leqslant j \leqslant L} \max_{0 \leqslant i < K} a_i^{(j)}(x_d) Z_{t-i}^{(j)}+\varepsilon_t(x_d)\qquad (1 \leqslant d \leqslant D)
\end{equation}
where
$(Z_t)_{t\in\Z}$ are independent unit-Fr\'echet,
$\varepsilon_t(x_d)$ are independent $N(0,\sigma^2)$ random variables,
the $a_i^{(j)}(x)$ are positive continuous functions defined on $[0,1]^q$ and
for every $x\in[0,1]^q$ we have that $\sum_{j\in\Z}\sum_{i\geqslant0} a_i^{(j)}(x)=1$.

It is important to note that the profiles $(a_0^{(j)},\ldots,a_{K-1}^{(j)})$ contain not only the information about the shapes of the profile but also, according to \ref{relfreq}, their probability of occurrence. The shapes will be denoted $a_{i,0}^{(j)}$ and their frequencies of occurrence $f^{(j)}$. That is
\begin{equation} \label{pbfreq}
a_i^{(j)} = \alpha^{(j)}a_{i,0}^{(j)}
\end{equation}
where the coefficients $\alpha^{(j)}$ must be chosen so that $f^{(j)}=p^{(j)}$ in \ref{relfreq}.

The first step of the procedure is to estimate the length of the tail dependence $K$.
This is done in \ref{estK} taking the average size of the clusters of exceedances over a threshold.
Next blocks of extremes are selected to estimate number of patterns $L$, the shapes $a_{i,0}^{(j)}$ of the parameter
functions and their frequencies $f^{(j)}$. The algorithm to locate the blocks of extremes explained in
\ref{ChooseExt} is based on a multivariate approach. The value of $L$ is determined in \ref{estL} as the number of clusters among the chosen blocks of extremes. The functions $a_{i,0}^{(j)}$ are yielded by the natural output (centroids, medoids,~...) of the partitioning algorithm used to determine $L$ and the values $f^{(j)}$ through the size of the different clusters. Then the solution $a_{i}^{(j)}$ of \eqref{pbfreq} is obtained in \ref{estA} thanks to an iterative algorithm. To measure the quality of the estimation, the quantification of the dissimilarity between the the original parameter functions and their estimations is done in \ref{hDist} in terms of the Hausdorff distance.

\subsection{Estimation of the length of the tail dependence ($K$)} \label{estK}

Eight estimators of $K$ have been tested for the model \eqref{practicalModel} with $\sigma=0$
when the knowledge of the parameter functions $a_i^{(j)}$ is replaced by the range $C$ given in \eqref{bounds}.

The first step is to select the values considered as extremes. A first approach consist of working
on $(\|X_t\|_\infty)$ and choosing the values above the threshold $\max_t(\|X_t\|_\infty)/C$. This will be referred
as the scalar version. A second approach can be to use all the available information
by doing the previous operation for the $D$ components of $(X_t)$ separately. In this last case
the threshold also depends on the location. This method will be referred as the multivariate version.

Once the extremes are selected, a runs declustering generates a sequence $(s_n)$ of all sizes of clusters of
extremes found in the univariate or multivariate scan. More precisely, only contiguous extremes were considered here
to make a cluster. This is the runs declustering with $r=0$.

From the sequence $(s_n)$, we estimate $K$ through $\operatorname{mean}(s_n)$, $\operatorname{median}(s_n)$
or $\operatorname{mode}(s_n)$ with the nomenclature as follow.

\begin{center}
\begin{tabular}{|l|l|l|l|}
\hline
& Average & Time series & Threshold (scalar or vector)\\
\hline
$\hat{K}_\mu$ & mean & $(\|X_t\|_\infty)$ & $\max_{t}((\|X_t\|_\infty)/C$\\
\hline
$\hat{K}_\mu^M$ & mean &  $(X_t)$ & $ (\max_{t}(X_t(d))/C)_{1 \leqslant d \leqslant D}$\\
\hline
$\hat{K}_m$ & median &  $(\|X_t\|_\infty)$ & $\max_{t}(\|X_t\|_\infty)/C$\\
\hline
$\hat{K}_m^M$ & median & $(X_t)$ & $(\max_{t}(X_t(d))/C)_{1 \leqslant d \leqslant D}$\\
\hline
$\hat{K}_o$ & mode &  $(\|X_t\|_\infty)$ & $\max_{t}(\|X_t\|_\infty)/C$\\
\hline
$\hat{K}_o^M$ & mode & $(X_t)$ & $(\max_{t}(X_t(d))/C)_{1 \leqslant d \leqslant D}$\\
\hline
\end{tabular}
\end{center}

The $\operatorname{ceil}$ or $\operatorname{floor}$ options to get rid of the decimals are used to build the eight following estimators:
$$ \begin{array}{l@{}l@{}l@{}l}
\hat{K}_1 = \operatorname{ceil}(K_\mu), &\hat{K}_3 = \operatorname{ceil}({K_\mu^M}), & \hat{K}_5 = \operatorname{ceil}(K_m),\ & \hat{K}_7 = {\hat{K}_o},\\

\hat{K}_2 = \operatorname{round}({K_\mu}),\ &\hat{K}_4 = \operatorname{round}({K_\mu^M}),\ &\hat{K}_6 = \operatorname{ceil}({K_m^M}),\ & \hat{K}_8 = {\hat{K}_o^M}.\\

\end{array}$$

\begin{figure}[here]
\begin{center}
\begin{pspicture}(-.5,-.5)(10.5,4.5)
\psset{xunit=1cm,yunit=4cm,runit=1cm}
\psline[linewidth=.5pt,plotstyle=curve,linecolor=lightgray]{-}(0,0.1)(10,0.1)
\psline[linewidth=.5pt,plotstyle=curve,linecolor=lightgray]{-}(0,0.2)(10,0.2)
\psline[linewidth=.5pt,plotstyle=curve,linecolor=lightgray]{-}(0,0.3)(10,0.3)
\psline[linewidth=.5pt,plotstyle=curve,linecolor=lightgray]{-}(0,0.4)(10,0.4)
\psline[linewidth=.5pt,plotstyle=curve,linecolor=lightgray]{-}(0,0.5)(10,0.5)
\psline[linewidth=.5pt,plotstyle=curve,linecolor=lightgray]{-}(0,0.6)(10,0.6)
\psline[linewidth=.5pt,plotstyle=curve,linecolor=lightgray]{-}(0,0.7)(10,0.7)
\psline[linewidth=.5pt,plotstyle=curve,linecolor=lightgray]{-}(0,0.8)(10,0.8)
\psline[linewidth=.5pt,plotstyle=curve,linecolor=lightgray]{-}(0,0.9)(10,0.9)
\psline[linewidth=.5pt,plotstyle=curve,linecolor=lightgray]{-}(0,1)(10,1)
\psline[linewidth=.5pt,plotstyle=curve,linecolor=lightgray]{-}(1,0)(1,1)
\psline[linewidth=.5pt,plotstyle=curve,linecolor=lightgray]{-}(2,0)(2,1)
\psline[linewidth=.5pt,plotstyle=curve,linecolor=lightgray]{-}(3,0)(3,1)
\psline[linewidth=.5pt,plotstyle=curve,linecolor=lightgray]{-}(4,0)(4,1)
\psline[linewidth=.5pt,plotstyle=curve,linecolor=lightgray]{-}(5,0)(5,1)
\psline[linewidth=.5pt,plotstyle=curve,linecolor=lightgray]{-}(6,0)(6,1)
\psline[linewidth=.5pt,plotstyle=curve,linecolor=lightgray]{-}(7,0)(7,1)
\psline[linewidth=.5pt,plotstyle=curve,linecolor=lightgray]{-}(8,0)(8,1)
\psline[linewidth=.5pt,plotstyle=curve,linecolor=lightgray]{-}(9,0)(9,1)
\psline[linewidth=.5pt,plotstyle=curve,linecolor=lightgray]{-}(10,0)(10,1)
\psline[linewidth=1pt,plotstyle=curve,linecolor=black]{->}(0,0)(10.5,0)
\psline[linewidth=1pt,plotstyle=curve,linecolor=black]{->}(0,0)(0,1.1)
\psline[linewidth=1pt,plotstyle=curve,linecolor=red,linestyle=dashed]{-}(0,.37408)(1,.49162)(2,.6509)(3,.73886)(4,.80132)(5,.79912)(6,.79782)(7,.79854)(8,.79374)(9,.7972)(10,.79752) 
\psline[linewidth=1pt,plotstyle=curve,linecolor=red]{-}(0,.36112)(1,.4543)(2,.58556)(3,.67616)(4,.81426)(5,.84202)(6,.85038)(7,.85518)(8,.86112)(9,.86508)(10,.86588) 
\psline[linewidth=1pt,plotstyle=curve,linecolor=blue,linestyle=dashed]{-}(0,.39884)(1,.51514)(2,.61628)(3,.61248)(4,.5861)(5,.58378)(6,.58142)(7,.58222)(8,.57872)(9,.58094)(10,.58144)  
\psline[linewidth=1pt,plotstyle=curve,linecolor=blue]{-}(0,.35648)(1,.4613)(2,.58852)(3,.63264)(4,.6586)(5,.66814)(6,.66886)(7,.67498)(8,.67772)(9,.68076)(10,.681)  
\psline[linewidth=1pt,plotstyle=curve,linecolor=green]{-}(0,.38292)(1,.5324)(2,.73942)(3,.83938)(4,.9002)(5,.9037)(6,.9037)(7,.90596)(8,.9057)(9,.90688)(10,.907) 
\psline[linewidth=1pt,plotstyle=curve,linecolor=magenta]{-}(0,.42718)(1,.5782)(2,.68938)(3,.68464)(4,.67018)(5,.6705)(6,.66488)(7,.6671)(8,.66658)(9,.66962)(10,.66756) 
\psline[linewidth=1pt,plotstyle=curve,linecolor=cyan]{-}(0,.3802)(1,.52942)(2,.76048)(3,.87274)(4,.9455)(5,.95158)(6,.95536)(7,.95614)(8,.95696)(9,.95882)(10,.95828) 
\psline[linewidth=1pt,plotstyle=curve,linecolor=yellow]{-}(0,.41374)(1,.56676)(2,.69512)(3,.72794)(4,.75384)(5,.7622)(6,.75646)(7,.76142)(8,.75972)(9,.76476)(10,.761) 
\psline[linewidth=1pt,plotstyle=curve,linecolor=red]{-}(1,.25)(2,.25) 
\psline[linewidth=1pt,plotstyle=curve,linecolor=red,linestyle=dashed]{-}(1,.15)(2,.15) 
\psline[linewidth=1pt,plotstyle=curve,linecolor=blue]{-}(3,.25)(4,.25) 
\psline[linewidth=1pt,plotstyle=curve,linecolor=blue,linestyle=dashed]{-}(3,.15)(4,.15) 
\psline[linewidth=1pt,plotstyle=curve,linecolor=green]{-}(5,.25)(6,.25) 
\psline[linewidth=1pt,plotstyle=curve,linecolor=yellow]{-}(7,.15)(8,.15) 
\psline[linewidth=1pt,plotstyle=curve,linecolor=cyan]{-}(7,.25)(8,.25) 
\psline[linewidth=1pt,plotstyle=curve,linecolor=magenta]{-}(5,.15)(6,.15) 
\uput[d](0,0){\footnotesize{$10$}}
\uput[d](1,0){\footnotesize{$20$}}
\uput[d](2,0){\footnotesize{$50$}}
\uput[d](3,0){\footnotesize{$100$}}
\uput[d](4,0){\footnotesize{$500$}}
\uput[d](5,0){\footnotesize{$1000$}}
\uput[d](6,0){\footnotesize{$1500$}}
\uput[d](7,0){\footnotesize{$2500$}}
\uput[d](8,0){\footnotesize{$5000$}}
\uput[d](9,0){\footnotesize{$7500$}}
\uput[d](10,0){\footnotesize{$10000$}}
\uput[l](0,.1){\footnotesize{$0.10$}}
\uput[l](0,.2){\footnotesize{$0.20$}}
\uput[l](0,.3){\footnotesize{$0.30$}}
\uput[l](0,.4){\footnotesize{$0.40$}}
\uput[l](0,.5){\footnotesize{$0.50$}}
\uput[l](0,.6){\footnotesize{$0.60$}}
\uput[l](0,.7){\footnotesize{$0.70$}}
\uput[l](0,.8){\footnotesize{$0.80$}}
\uput[l](0,.9){\footnotesize{$0.90$}}
\uput[l](0,1.00){\footnotesize{$1.00$}}
\uput[dr](0.1,1.1){\footnotesize{$\hat{p}$}}
\uput[u](10.4,0){\footnotesize{$T$}}
\uput[r](1.9,.25){\footnotesize{$\hat{K}_1$}}
\uput[r](1.9,.15){\footnotesize{$\hat{K}_2$}}
\uput[r](3.9,.25){\footnotesize{$\hat{K}_3$}}
\uput[r](3.9,.15){\footnotesize{$\hat{K}_4$}}

\uput[r](5.9,.25){\footnotesize{$\hat{K}_5$}}
\uput[r](7.9,.15){\footnotesize{$\hat{K}_8$}}
\uput[r](7.9,.25){\footnotesize{$\hat{K}_7$}}

\uput[r](5.9,.15){\footnotesize{$\hat{K}_6$}}
\end{pspicture}
\caption{Proportion of success in estimating $K$.}
\label{ComparisonK}
\end{center}
\end{figure}
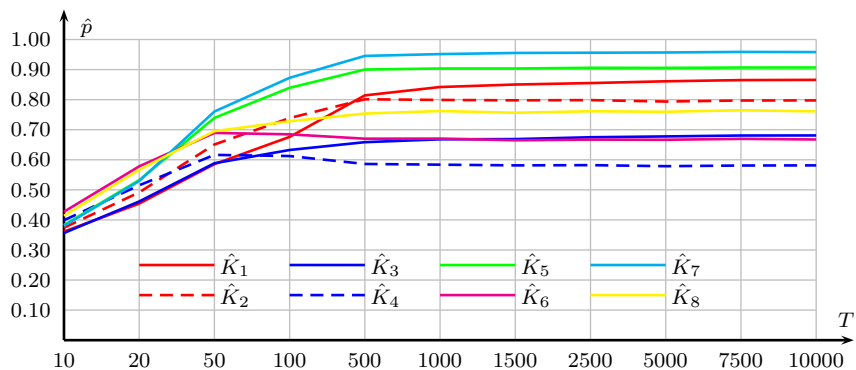

Figure \ref{ComparisonK} shows the success rate the eight estimators of $K$
against the length of the simulated chain. The tests were performed with $N=50000$ trials at each step: for $C$ from 1 to 10, $D$ from 1 to 20, $K$ from 1 to 5,  $L$ from 1 to 5, $\sigma=0$ and, for each of these parameters, 10 different sets of coefficients $a_i^{(j)}$ randomly generated (uniformly, without time or space correlation).

According to these empirical results, the winner for $T\geqslant35$ is the univariate version and the mode as average cluster size. If $T\leqslant35$ the best success rate is obtained with the multivariate version and the median.

\subsection{Extremal clustering} \label{ChooseExt}

Once we know the length $K$ of the tail dependence thanks to \ref{estK}, the next step of the procedure studied here to recover the parameter functions $a_i^{(j)}$ of a theoretical CM3 process \eqref{theoreticalModel} is to locate the blocks of extremes. Indeed, according to \ref{correlation}, the probability that at least one block of length $K$ in the chain is block profile of type $l$ in a sample of size $T$ is greater than
$$ 1-(1-p^{(l)})^{T-K+1} $$
which tends to 1 as $T$ tends to infinity.

The suggested method locates the positions of blocks of extremes in the practical model \eqref{practicalModel} maximizing the ``likelihood'' of being a multivariate extreme. This idea comes from the wish not to lose information across the $D$ dimensions. Nevertheless a bad situation can still happen when, for some $0 \leqslant i^\ast < K$, all $a_{i^\ast}^{(j)}(x_d)$ are negligible in comparison with the $a_{i}^{(j)}(x_d)$, $i\neq i^\ast$, for instance. If the points $x_1,\ldots,x_D$ of the discretization are far enough to obtain independent-like patterns, it is unlikely that all the $a_{i^\ast}^{(j)}(x_d)$ are negligible in the same time.

We explain the method on the following example with $D=2$ and $K=3$:
$$\begin{array}{|c|ccccccc|}
\hline
X_t(1) & 5 & 3 & 4 & 14 & 19 & 2 & 7 \\
\hline
X_t(2) & 6 & 1 & 10 & 5 & 2 & 1 & 4\\
\hline
\end{array}.$$

\noindent\underline{First step}

\nopagebreak
\vspace{.09cm}
\nopagebreak
\noindent  Using the order statistics, mark the $K$ largest values in the $D$ chains by 1.
$$\begin{array}{|c|ccccccc|}
\hline
d=1 & 0 & 0 & 0 & 1 & 1 & 0 & 1 \\
\hline
d=2 & 1 & 0 & 1 & 1 & 0 & 0 & 0\\
\hline
\end{array}$$

\noindent\underline{Second step}

\nopagebreak
\vspace{.09cm}
\nopagebreak
\noindent  Compute the sum of the extremal status for each $t$.
$$\begin{array}{|c|ccccccc|}
\hline
d=1 & 0 & 0 & 0 & 1 & 1 & 0 & 1 \\
\hline
d=2 & 1 & 0 & 1 & 1 & 0 & 0 & 0\\
\hline
\Rightarrow\lambda    & 1 & 0 & 1 & 2 & 1 & 0 & 1\\
\hline
\end{array}$$

\noindent\underline{Third step}

\nopagebreak
\vspace{.09cm}
\nopagebreak
\noindent  Compute the moving sum (MS) of order $K$. This is considered as the likelihood $\lambda$ to have an large value at time $t$ among the $(Z_{t}^{(l)})_{1 \leqslant l \leqslant L}$.
$$\begin{array}{|c|ccccccc|}
\hline
\lambda    & 1 & 0 & 1 & 2 & 1 & 0 & 1\\
\hline
MS   & 1 & 1 & 2 & 3 & 4 & 3 & 2\\
\hline
\end{array}$$

\noindent\underline{Extract the profile}

\nopagebreak
\vspace{.09cm}
\nopagebreak
\noindent  Find the index $t$ that maximizes the moving sum. Then
\begin{equation}\label{sto}
(X_{t-K+1},\ldots,X_t)/\|X_{t-K+1}\|_\infty
\end{equation}
is the shape of the first block profile to store in the memory:

$$\begin{array}{|c|ccccccc|}
\hline
MS   & 1 & 1 & 2 & 3 & \red{\textbf{4}} & 3 & 2\\
\hline
X_t(1) & 5 & 3 & \red{4} & \red{14} & \red{19} & 2 & 7 \\
\hline
X_t(2) & 6 & 1 & \red{\underline{10}} & \red{5} & \red{2} & 1 & 4\\
\hline
\end{array}$$

$$\begin{array}{rrcccl}

\textrm{Shape}_1(x_1)=&(&4/\underline{10},& 14/\underline{10},& 19/\underline{10}&)\\

\textrm{Shape}_1(x_2)=&(&10/\underline{10},& 5/\underline{10},& 2/\underline{10}&)\\

\end{array}$$
To decide between multiple maxima, for instance if
$$\begin{array}{|c|ccccccc|}
\hline
MS   & 1 & 1 & 2 & \red{\textbf{4}} & 2 & \red{\textbf{4}} & \red{\textbf{4}}\\
\hline
\end{array}$$
we first choose the single maximum. If there are consecutive maxima, as a second criterion, we take the block that maximizes $\sum_{t\in\textrm{block}}\|X_t\|_\infty$ among those.\\

\noindent  Repeat this loop until having gathered the desired number $Q$ of time-disjoint blocks of extremes of length $K$ ($Q=\underline{p}T$ is the suggestion of \ref{relfreq} if we only know a uniform bound $C$ on the variation of the parameter functions).

\begin{center}
\begin{figure}[here!]
\includegraphics[scale=.5]{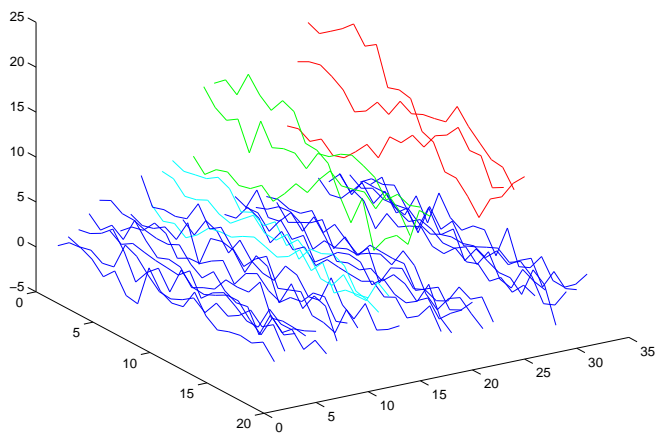}
\caption{Blocks of extremes obtained for a CM3 with measurement errors ($\sigma=1$).}
\end{figure}
\end{center}

\subsection{Estimation of the number of patterns ($L$)} \label{estL}

With $Q$ blocks of extremes of length $K$ normalized as in \eqref{sto}, the goal is to estimate the number of reproducible patterns $L$ in the observed process. To do this, we create a $Q\times KD$-table inside which each of the $Q$ lines is made of the $D$ temporal vectors of length $K$ placed successively. We estimate $L$ with the number of clusters for the observations of the table.\\

\noindent\underline{Partitioning methods}

\nopagebreak
\vspace{.09cm}
\nopagebreak
To break the lines of the table up into groups, we tried several algorithms among which five retained our attention: hierarchical clustering with Ward's aggregation criterion, hierarchical clustering with the Euclidean distance between the centroids \cite{HTR09,WAR63}, $k$-means with the Euclidean squared distance, $k$-means with Pearson's correlation after standardization \cite{SEB84,SPA85} and finally Partitioning Around Medoids (PAM) with the classical Euclidean distance  \cite{TK06,VPB03}.\\

\noindent\underline{Number of clusters}

\begin{figure}[here!]
\begin{center}
\includegraphics[scale=0.4]{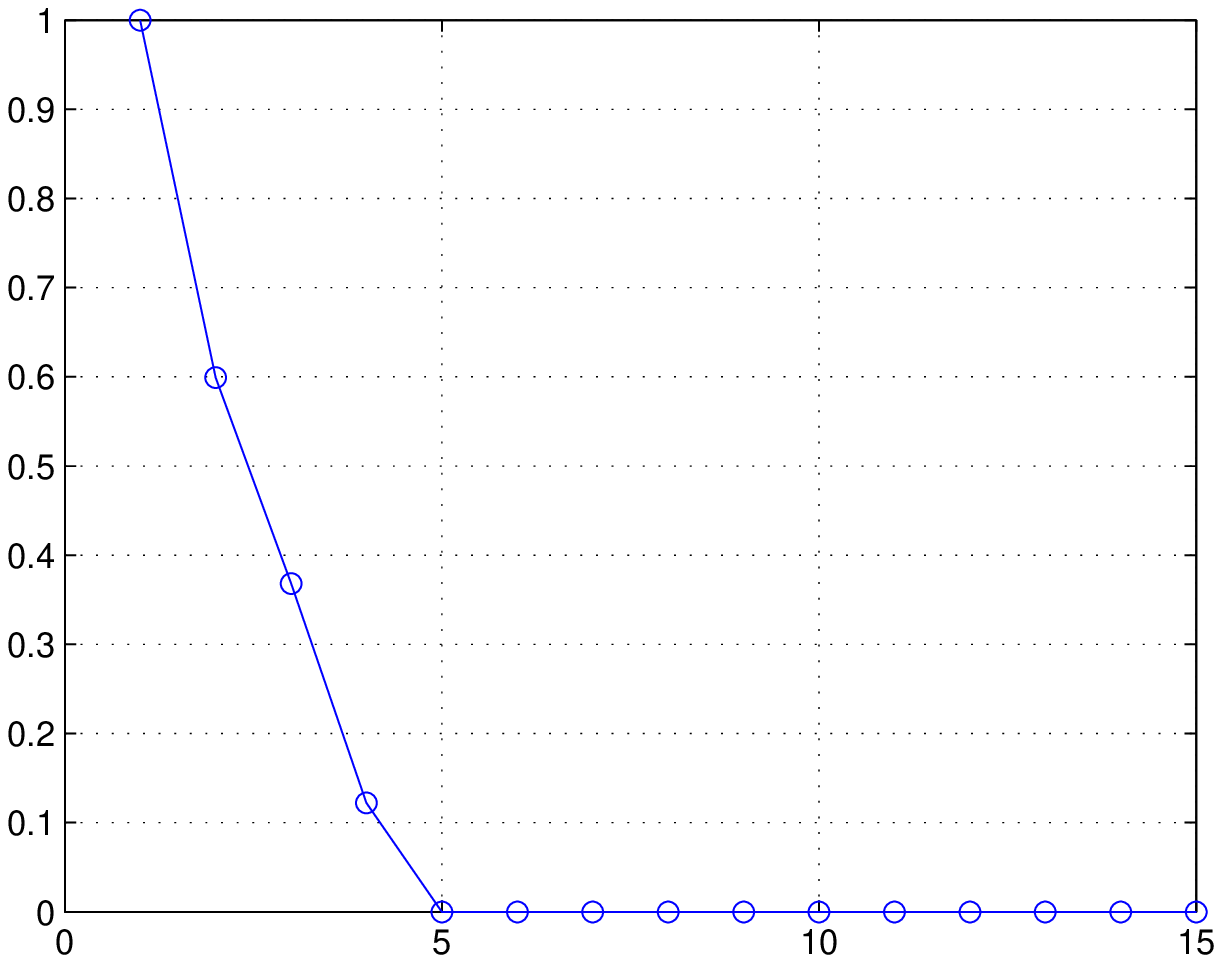}
\includegraphics[scale=0.4]{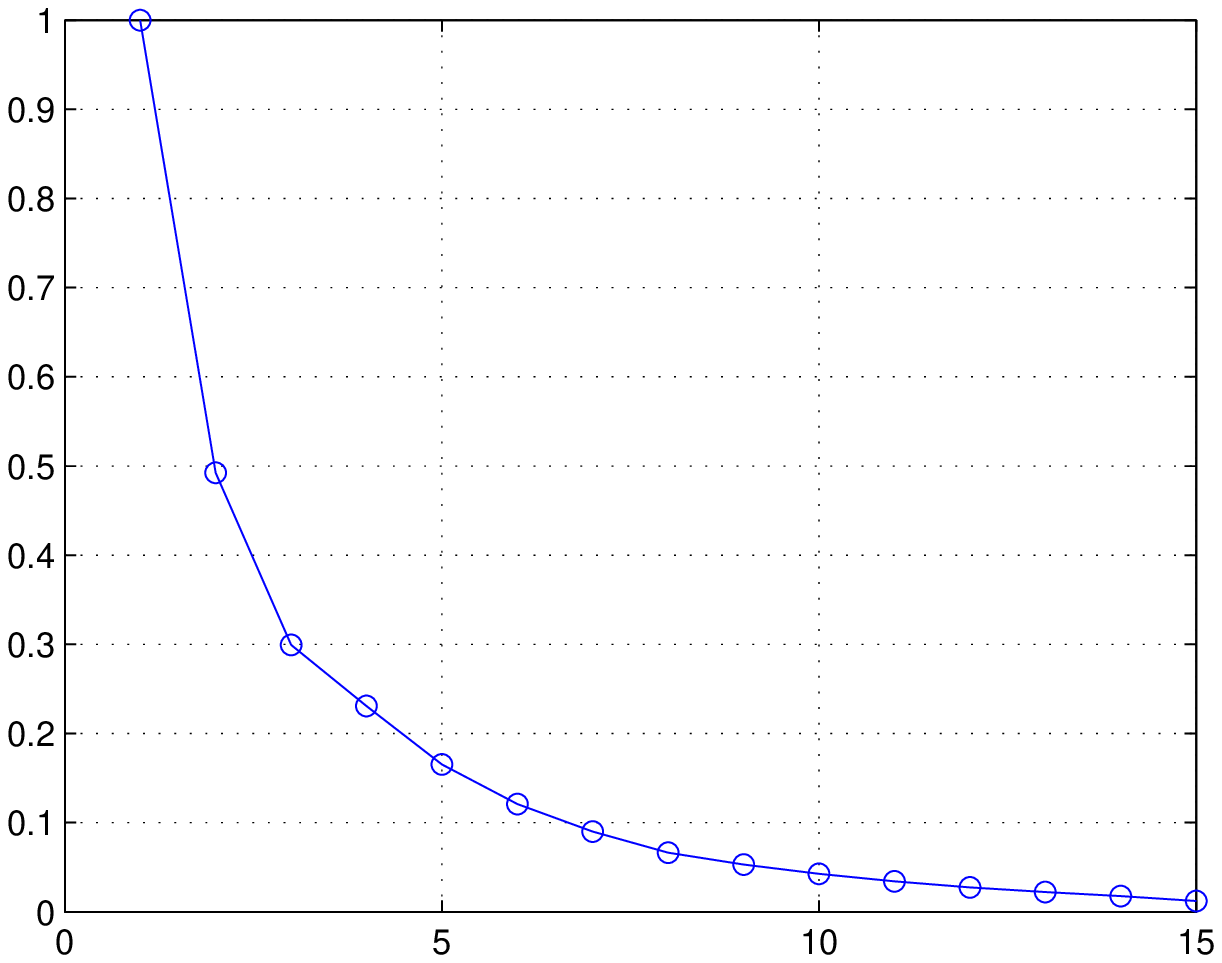}
\begin{pspicture}(-7,0)(7,0)
\psset{xunit=1cm,yunit=1cm,runit=1cm}
\uput[u](-3.8,0){\footnotesize{Number of clusters}}
\uput[u](2.2,0){\footnotesize{Number of clusters}}
\uput[u](-6.75,3.9){\footnotesize{$\frac{\textrm{SSE}}{\textrm{SStot}}$}}
\uput[u](-.75,3.9){\footnotesize{$\frac{\textrm{SSE}}{\textrm{SStot}}$}}
\psline[linewidth=.5pt,linecolor=red]{-}(-6.1,1.549)(-1.51,1.549)
\psline[linewidth=.5pt,linecolor=red]{-}(-.085,1.549)(4.51,1.549)
\end{pspicture}
\caption{Elbow: With $L=5$, perfect clustering on the left, CM3 with errors on the right.}
\label{FigEl}
\end{center}
\end{figure}

\nopagebreak
\vspace{.09cm}
\nopagebreak
For each of those algorithms we implemented two criteria to determine the number of clusters. The first: one
stops when the percentage of the total variance not explained by the clustering
is less than $20\%$, i.e.~when
$$\frac{\textrm{SSE}}{\textrm{SStot}}= \frac{\sum_{\operatorname{cluster}}(\operatorname{nbobs}(\operatorname{cluster})-1) \sum_{\textrm{variable in cluster}}s^2(\operatorname{variable})}{(Q -1)\sum_{\textrm{variable in table}}s^2(\operatorname{variable})}\leqslant0.20. $$
We refer to this method as the \emph{elbow method} \cite{KS96,MAR79} (see Figure~\ref{FigEl}).

\begin{figure}[here!]
\begin{center}
\includegraphics[scale=0.4]{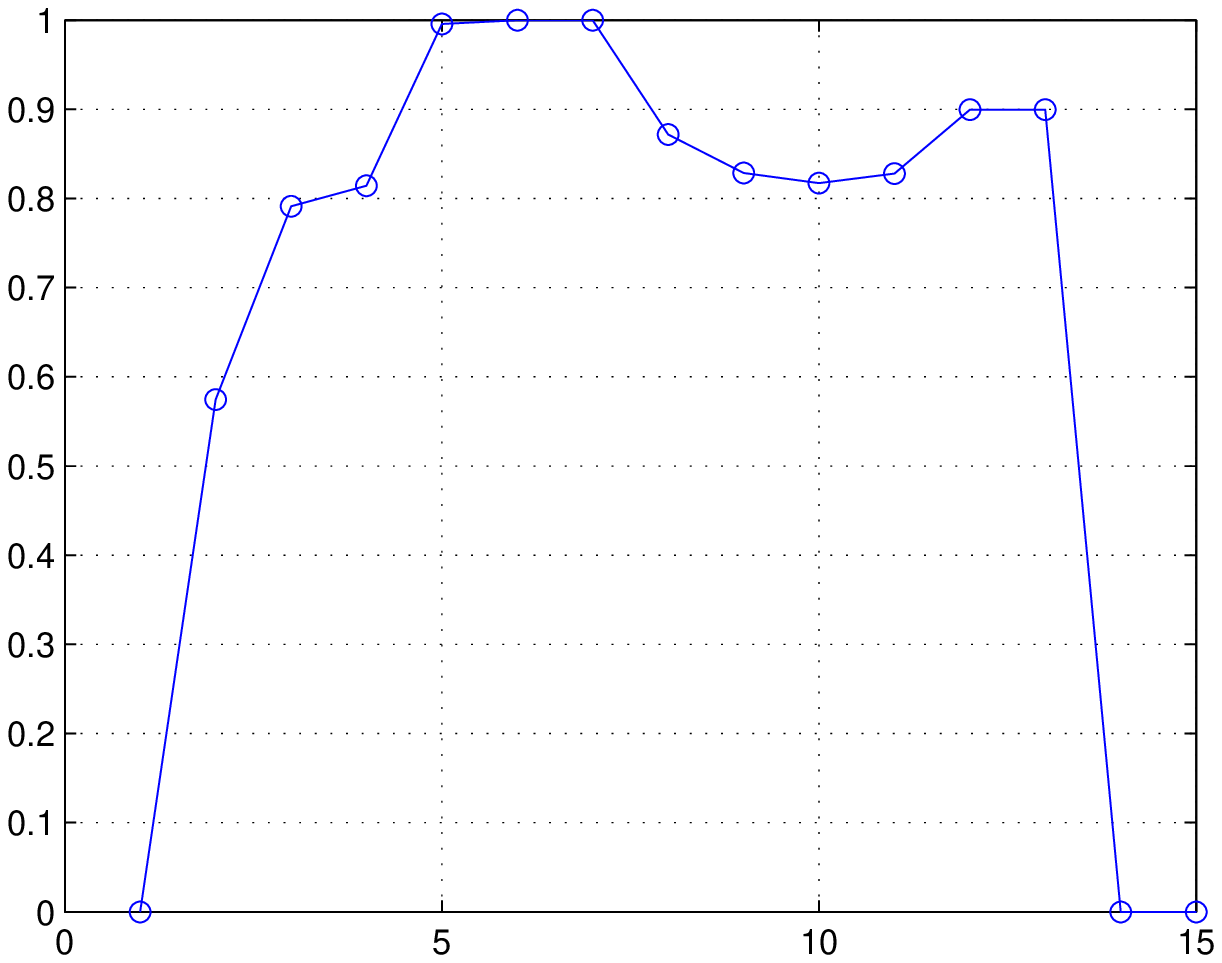}
\includegraphics[scale=0.4]{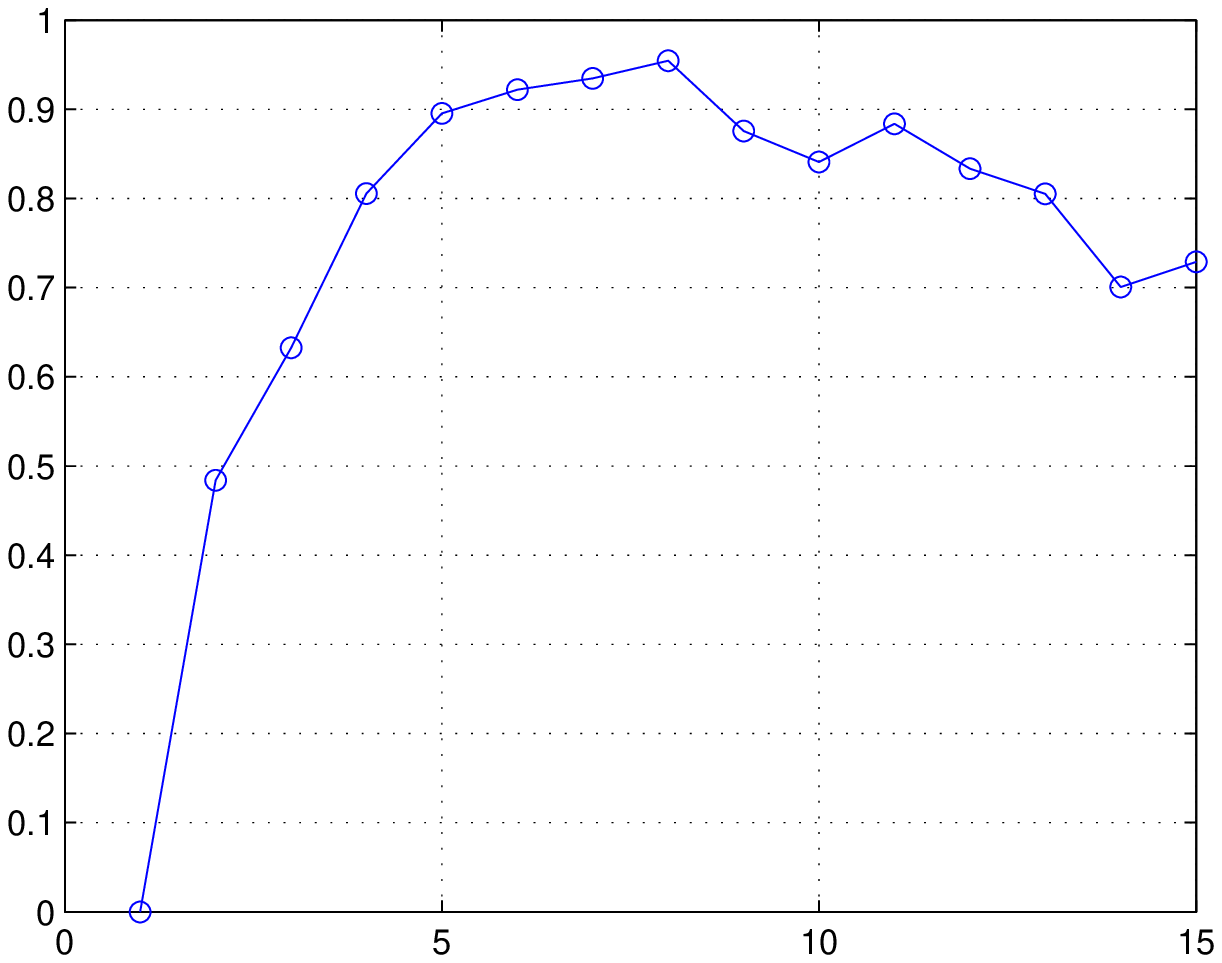}
\begin{pspicture}(-7,0)(7,0)
\psset{xunit=1cm,yunit=1cm,runit=1cm}
\uput[u](-3.8,0){\footnotesize{Number of clusters}}
\uput[u](2.2,0){\footnotesize{Number of clusters}}
\uput[u](-6.75,3.9){\footnotesize{$\frac{\textrm{TtSil}}{Q}$}}
\uput[u](-.75,3.9){\footnotesize{$\frac{\textrm{TtSil}}{Q}$}}
\psline[linewidth=.5pt,linecolor=green]{-}(-6.1,3.9)(-1.51,3.9)
\psline[linewidth=.5pt,linecolor=green]{-}(-.085,3.9)(4.51,3.9)
\end{pspicture}
\caption{Silhouette: Perfect clustering on the left, CM3 with errors on the right ($L=5$).}
\label{FigSil}
\end{center}
\end{figure}

The second method to find the number of clusters here is the first value that
yields a total silhouette TtSil for the clustering above $85\%$ of $Q$.
Let $a(q)$ be the average
distance between the $q^{\textrm{th}}$ observation and the members of its own
cluster. Then repeat this operation between
the $q^{\textrm{th}}$ observation and the members of all the other clusters,
and set $b(q)$ to the lowest value found.
The silhouette $s(q)$ of the $q^{\textrm{th}}$ observation is $$
s(q)=\frac{b(q)-a(q)}{\max\{a(q),b(q)\}}. $$ 
Thus $-1\leqslant s(q)\leqslant1$ and $s(q)$ measures how dissimilar the
$q^{\textrm{th}}$ observation is to its own cluster \cite{KR90,LO04,ROU87}.
The distance taken into account here is the Euclidean squared distance.
We stop the partitioning
at the smallest number of clusters satisfying
$$ \frac{\textrm{TtSil}}{Q} = \frac{1}{Q} \sum_{q=1}^Q s(q)\geqslant 0.85 $$
if this occurs. We refer to this method as the \emph{silhouette method}
(see Figure~\ref{FigSil}).

Both methods are unable to detect that $L=1$. We have thus considered $\hat{L}=1$ when the estimated variance of each variable is less than 0.005.\\

\noindent\underline{Estimation}

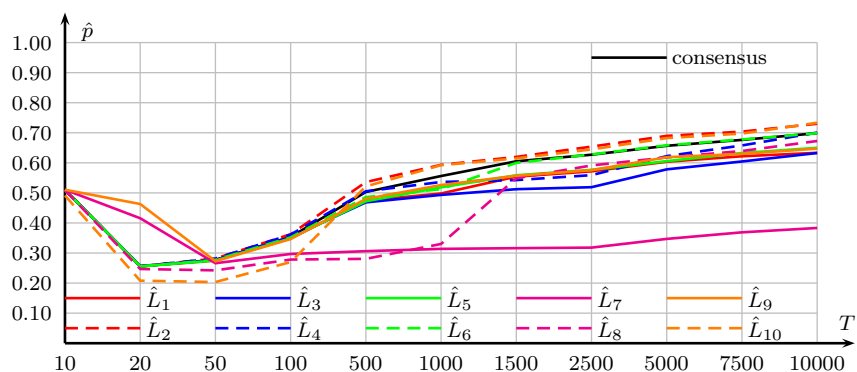
\begin{figure}[here!]
\begin{center}
\begin{pspicture}(-.5,-.2)(10.5,4.5)
\psset{xunit=1cm,yunit=4cm,runit=1cm}
\psline[linewidth=.5pt,plotstyle=curve,linecolor=lightgray]{-}(0,0.1)(10,0.1)
\psline[linewidth=.5pt,plotstyle=curve,linecolor=lightgray]{-}(0,0.2)(10,0.2)
\psline[linewidth=.5pt,plotstyle=curve,linecolor=lightgray]{-}(0,0.3)(10,0.3)
\psline[linewidth=.5pt,plotstyle=curve,linecolor=lightgray]{-}(0,0.4)(10,0.4)
\psline[linewidth=.5pt,plotstyle=curve,linecolor=lightgray]{-}(0,0.5)(10,0.5)
\psline[linewidth=.5pt,plotstyle=curve,linecolor=lightgray]{-}(0,0.6)(10,0.6)
\psline[linewidth=.5pt,plotstyle=curve,linecolor=lightgray]{-}(0,0.7)(10,0.7)
\psline[linewidth=.5pt,plotstyle=curve,linecolor=lightgray]{-}(0,0.8)(10,0.8)
\psline[linewidth=.5pt,plotstyle=curve,linecolor=lightgray]{-}(0,0.9)(10,0.9)
\psline[linewidth=.5pt,plotstyle=curve,linecolor=lightgray]{-}(0,1)(10,1)
\psline[linewidth=.5pt,plotstyle=curve,linecolor=lightgray]{-}(1,0)(1,1)
\psline[linewidth=.5pt,plotstyle=curve,linecolor=lightgray]{-}(2,0)(2,1)
\psline[linewidth=.5pt,plotstyle=curve,linecolor=lightgray]{-}(3,0)(3,1)
\psline[linewidth=.5pt,plotstyle=curve,linecolor=lightgray]{-}(4,0)(4,1)
\psline[linewidth=.5pt,plotstyle=curve,linecolor=lightgray]{-}(5,0)(5,1)
\psline[linewidth=.5pt,plotstyle=curve,linecolor=lightgray]{-}(6,0)(6,1)
\psline[linewidth=.5pt,plotstyle=curve,linecolor=lightgray]{-}(7,0)(7,1)
\psline[linewidth=.5pt,plotstyle=curve,linecolor=lightgray]{-}(8,0)(8,1)
\psline[linewidth=.5pt,plotstyle=curve,linecolor=lightgray]{-}(9,0)(9,1)
\psline[linewidth=.5pt,plotstyle=curve,linecolor=lightgray]{-}(10,0)(10,1)
\psline[linewidth=1pt,plotstyle=curve,linecolor=black]{->}(0,0)(10.5,0)
\psline[linewidth=1pt,plotstyle=curve,linecolor=black]{->}(0,0)(0,1.1)
\psline[linewidth=1pt,plotstyle=curve,linecolor=black]{-}(0,.5105)(1,.25616)(2,.277)(3,.35367)(4,.50317)(5,.55617)(6,.60517)(7,.627)(8,.65667)(9,.67633)(10,.69867)
\psline[linewidth=1pt,plotstyle=curve,linecolor=red]{-}(0,.5105)(1,.25584)(2,.27517)(3,.34867)(4,.47417)(5,.498)(6,0.55417)(7,.57167)(8,.60367)(9,.62183)(10,.63333)
\psline[linewidth=1pt,plotstyle=curve,linecolor=red,linestyle=dashed]{-}(0,.5105)(1,.256)(2,.28167)(3,.3625)(4,.536)(5,.59333)(6,.61967)(7,.654)(8,.68967)(9,.70383)(10,.73)
\psline[linewidth=1pt,plotstyle=curve,linecolor=blue]{-}(0,.5105)(1,.256)(2,.27717)(3,.35283)(4,.4685)(5,.49333)(6,.51233)(7,.519)(8,.57833)(9,.60450)(10,.63267)
\psline[linewidth=1pt,plotstyle=curve,linecolor=blue,linestyle=dashed]{-}(0,.5105)(1,.256)(2,.28167)(3,.36133)(4,.50533)(5,.53633)(6,.54267)(7,.559)(8,.62233)(9,.65783)(10,.70067)
\psline[linewidth=1pt,plotstyle=curve,linecolor=green]{-}(0,.5105)(1,.25584)(2,.275)(3,.34883)(4,.47217)(5,.5185)(6,.55867)(7,.57667)(8,.60567)(9,.63267)(10,.65)
\psline[linewidth=1pt,plotstyle=curve,linecolor=green,linestyle=dashed]{-}(0,.5105)(1,.25506)(2,.27733)(3,.352)(4,.486)(5,.51267)(6,.59983)(7,.628)(8,.65833)(9,.67783)(10,.69867)
\psline[linewidth=1pt,plotstyle=curve,linecolor=magenta]{-}(0,.5105)(1,.41569)(2,.26617)(3,.29683)(4,.30617)(5,.31383)(6,.31617)(7,.318)(8,.347)(9,.36883)(10,.38333)
\psline[linewidth=1pt,plotstyle=curve,linecolor=magenta,linestyle=dashed]{-}(0,.5095)(1,.2469)(2,.24233)(3,.27833)(4,.28067)(5,.33017)(6,.54867)(7,.59133)(8,.618)(9,.63933)(10,.67267)
\psline[linewidth=1pt,plotstyle=curve,linecolor=orange]{-}(0,.5105)(1,.46306)(2,.27283)(3,.3465)(4,.48)(5,.52617)(6,.557)(7,.57567)(8,.61867)(9,.62967)(10,.64733)
\psline[linewidth=1pt,plotstyle=curve,linecolor=orange,linestyle=dashed]{-}(0,.49133)(1,.208)(2,.2035)(3,.26983)(4,.52167)(5,.59217)(6,.61483)(7,.64467)(8,.68233)(9,.69733)(10,.73333)

\psline[linewidth=1pt,plotstyle=curve,linecolor=red]{-}(0,.15)(1,.15)
\psline[linewidth=1pt,plotstyle=curve,linecolor=red,linestyle=dashed]{-}(0,.05)(1,.05)
\psline[linewidth=1pt,plotstyle=curve,linecolor=blue]{-}(2,.15)(3,.15)
\psline[linewidth=1pt,plotstyle=curve,linecolor=blue,linestyle=dashed]{-}(2,.05)(3,.05)
\psline[linewidth=1pt,plotstyle=curve,linecolor=green]{-}(4,.15)(5,.15)
\psline[linewidth=1pt,plotstyle=curve,linecolor=green,linestyle=dashed]{-}(4,.05)(5,.05)
\psline[linewidth=1pt,plotstyle=curve,linecolor=magenta]{-}(6,.15)(7,.15)
\psline[linewidth=1pt,plotstyle=curve,linecolor=magenta,linestyle=dashed]{-}(6,.05)(7,.05)
\psline[linewidth=1pt,plotstyle=curve,linecolor=orange]{-}(8,.15)(9,.15)
\psline[linewidth=1pt,plotstyle=curve,linecolor=orange,linestyle=dashed]{-}(8,.05)(9,.05)
\psline[linewidth=1pt,plotstyle=curve,linecolor=black]{-}(7,.95)(8,.95)
\uput[d](0,0){\footnotesize{$10$}}
\uput[d](1,0){\footnotesize{$20$}}
\uput[d](2,0){\footnotesize{$50$}}
\uput[d](3,0){\footnotesize{$100$}}
\uput[d](4,0){\footnotesize{$500$}}
\uput[d](5,0){\footnotesize{$1000$}}
\uput[d](6,0){\footnotesize{$1500$}}
\uput[d](7,0){\footnotesize{$2500$}}
\uput[d](8,0){\footnotesize{$5000$}}
\uput[d](9,0){\footnotesize{$7500$}}
\uput[d](10,0){\footnotesize{$10000$}}
\uput[l](0,.1){\footnotesize{$0.10$}}
\uput[l](0,.2){\footnotesize{$0.20$}}
\uput[l](0,.3){\footnotesize{$0.30$}}
\uput[l](0,.4){\footnotesize{$0.40$}}
\uput[l](0,.5){\footnotesize{$0.50$}}
\uput[l](0,.6){\footnotesize{$0.60$}}
\uput[l](0,.7){\footnotesize{$0.70$}}
\uput[l](0,.8){\footnotesize{$0.80$}}
\uput[l](0,.9){\footnotesize{$0.90$}}
\uput[l](0,1.00){\footnotesize{$1.00$}}
\uput[dr](0.1,1.1){\footnotesize{$\hat{p}$}}
\uput[u](10.4,0){\footnotesize{$T$}}
\uput[r](0.9,.15){\footnotesize{$\hat{L}_1$}}
\uput[r](0.9,.05){\footnotesize{$\hat{L}_2$}}
\uput[r](6.9,.15){\footnotesize{$\hat{L}_7$}}
\uput[r](2.9,.15){\footnotesize{$\hat{L}_3$}}
\uput[r](2.9,.05){\footnotesize{$\hat{L}_4$}}
\uput[r](8.9,.15){\footnotesize{$\hat{L}_9$}}
\uput[r](4.9,.15){\footnotesize{$\hat{L}_5$}}
\uput[r](4.9,.05){\footnotesize{$\hat{L}_{6}$}}
\uput[r](7.9,.945){\footnotesize{consensus}}
\uput[r](8.9,.05){\footnotesize{$\hat{L}_{10}$}}
\uput[r](6.9,.05){\footnotesize{$\hat{L}_8$}}
\end{pspicture}
\caption{Proportion of success in estimating $L$.}
\label{ComparisonL}
\end{center}
\end{figure}

\nopagebreak
\vspace{.09cm}
\nopagebreak
For $\sigma=0$, Figure \ref{ComparisonL} shows the success rate of the following eleven estimators of $L$ against the length of the simulated chain.
\begin{center}
\begin{tabular}{|l|l|l|l|}
\hline
& Clustering & Distance & Number of clusters\\
\hline
\hline
$\hat{L}_1$ & Hierarchical & Ward &  Elbow \\
\hline
$\hat{L}_2$ & Hierarchical & Ward & Silhouette \\
\hline
$\hat{L}_3$ & Hierarchical & Euclidean / centroids &  Elbow \\
\hline
$\hat{L}_4$ & Hierarchical & Euclidean / centroids & Silhouette \\
\hline
$\hat{L}_5$ & $k$-means & Euclidean squared & Elbow \\
\hline
$\hat{L}_6$ & $k$-means & Euclidean squared & Silhouette \\
\hline
$\hat{L}_7$ & $k$-means & Pearson's correlation &  Elbow \\
\hline
$\hat{L}_8$ & $k$-means & Pearson's correlation &  Silhouette \\
\hline
$\hat{L}_{9}$ & PAM & Euclidean &  Elbow \\
\hline
$\hat{L}_{10}$ & PAM & Euclidean &  Silhouette \\
\hline
$\hat{L}_{11}$ & \multicolumn{3}{l|}{Consensus: $\hat{L}_{11}=\operatorname{mode}(\hat{L}_1,\ldots,\hat{L}_{10})$}  \\
\hline
\end{tabular}
\end{center}
The tests were performed with
$N=6000$ trials at each step: for $C$ in $\{2,4,6,8,10\}$, $D$ in
$\{1,5,10,15,20\}$, $K$ from 2 to 5 and $L$ from 1 to 5. The number of blocks
of extremes chosen to build the table is
$Q=\min(\operatorname{ceil}(\underline{p}T),100)$ with $\underline{p}$ from
\ref{relfreq}. We added the constraint that it has to be possible
to see each profile once, i.e.~$(K+1)L\leqslant T$, to exclude
challenges such as finding 5 profiles of length 5 in a chain of length 10.

The combination between PAM and elbow yields the best success rate for small
sample sizes $T\leqslant 50$. For $T\geqslant 50$ Ward's algorithm with the
silhouette is the best. Here the consensus curve does not provide a better
performance.

The gap between $T=20$ and $T=500$ is an intermediate area between a situation
where the trivial case $Q=L$ frequently occurs and the apparition of asymptotic
properties.

These results are not so excellent but, since the implemented algorithm is the
``eye'' of the analyst, it only sees what we want it to see. For real applications
it probably does not matter if not very frequent profiles are missed. Moreover
optimizing $C$, $Q$ and the thresholds for a precise sample can
improve the performances for that particular situation.

\subsection{Recovering the parameter functions $a_i^{(j)}$} \label{estA}

Once wet got the length $K$ of the tail dependence in \ref{estK} and the number $L$ of different patterns in \ref{estL}, it remains to estimate the parameter functions $a_i^{(j)}$ ($0\leqslant i < K$, $1\leqslant j \leqslant L$). Depending on the partitioning method to estimate $L$, we choose for the shapes $a_{i,0}^{(j)}$ of the $L$ profiles the natural output of the algorithm: the centroids for the hierarchical clustering and $k$-means and the medoids with PAM. Essentially the difference is that in the first case the estimator of a shape is a mean of observations and in the second case a median observation.

The relationship between the parameter functions $a_i^{(j)}$, the shapes $a_{i,0}^{(j)}$ and their frequencies of occurrence is given by \eqref{profprob}. The theoretical probabilities of occurrence of the different profiles must match the empirical frequencies $(f^{(1)},\ldots,f^{(L)})$ in the table of \ref{estL}. Thus the last step is to normalize the $a_{i,0}^{(j)}$ to obtain $\hat{a}_i^{(j)}$ such that on the one hand
\begin{equation}\label{target1}
\harm\left((m_k^{(l;1)})_{k,l}\right)/f^{(1)}=\ldots=\harm\left((m_k^{(l;L)})_{k,l}\right)/f^{(L)}
\end{equation}
($1\leqslant k \leqslant 2K-1$, $1\leqslant l \leqslant L$)
and on the other hand
\begin{equation}\label{target2}
\sum_{j\in\Z}\sum_{i\geqslant0} \hat{a}_i^{(j)}(x_d)=1\qquad(1\leqslant d \leqslant D).
\end{equation}
We seek a relation of the type $\hat{a}_i^{(l)}=\alpha^{(l)} a_{i,0}^{(l)}$ ($1 \leqslant l \leqslant L $) which is a problem of rank $L$.
The solution can be obtained iteratively: if $a_{i,n}^{(j)}$ is the $n^{\textrm{th}}$ update of $a_{i,0}^{(j)}$ and $(p_n^{(1)},\ldots,p_n^{(L)})$ the probabilities given by \eqref{profprob} at step $n$, then define $a_{i,n+1}^{(j)}$ by
\begin{equation}\label{algo}
a_{i,n+1}^{(j)}(\cdot)=a_{i,n+1}^{(j)}(\cdot)\frac{f^{(j)}}{p_{n}^{(j)}}
\end{equation}
and stop when the error $|(f^{(1)},\ldots,f^{(L)})-(p_{n+1}^{(1)},\ldots,p_{n+1}^{(L)})|_\infty$ is small.
The iterative algorithm \eqref{algo}, if it converges, converges to the solution of \eqref{target1}. 
Indeed, let
$$ M_{n}^{(l^\star)}=\left(
\begin{array}{ccc}
m_{1,n}^{(1:l^\star)} & \ldots & m_{2K-1,n}^{(1:l^\star)}\\
\ldots&\ldots&\ldots\\
m_{1,n}^{(L:l^\star)} & \ldots & m_{2K-1,n}^{(L:l^\star)}\\
\end{array}
\right) $$
be the $n^{\textrm{th}}$ update of the initial collection of numbers $(m_{k,0}^{(l;l^\star)})_{k,l}$ given the $a_{i,0}^{(j)}$. Define
$$ T:\R^{(2K-1)L^2}\to\R^{(2K-1)L^2}:(M_n^{(1)},\ldots,M_n^{(L)})\mapsto (M_{n+1}^{(1)},\ldots,M_{n+1}^{(L)}) $$
as the resulting operation of \eqref{algo} on the $m_{k,n}^{(l;l^\star)}$. The operator $T$ acts so that $M_{n}^{(l^\star)}$ becomes
$$ M_{n+1}^{(l^\star)}=\left(
\begin{array}{ccl}
(m_{1,n}^{(1:l^\star)} & \ldots &\  m_{2K-1,n}^{(1:l^\star)})\frac{\displaystyle f^{(l^\star)} }{\displaystyle p_{n}^{(l^\star)}} \frac{\displaystyle p_{n}^{(1)} }{\displaystyle f^{(1)} }\\
\ldots&\ldots&\ \ldots\\
m_{1,n}^{(l^\star:l^\star)} & \ldots &\  m_{2K-1,n}^{(l^\star:l^\star)}\\
\ldots&\ldots&\ \ldots\\
(m_{1,n}^{(L:l^\star)} & \ldots &\  m_{2K-1,n}^{(L:l^\star)})\frac{\displaystyle f^{(l^\star)} }{\displaystyle p_{n}^{(l^\star)}} \frac{\displaystyle p_{n}^{(L)} }{\displaystyle f^{(L)} }\\
\end{array}
\right) .$$
Observe that \eqref{target1} is equivalent to $ T(M^{(1)},\ldots,M^{(L)})=(M^{(1)},\ldots,M^{(L)}) $. Consequently, if $(M_n^{(1)},\ldots,M_n^{(L)})_n$ converges, since $T$ is continuous, \eqref{algo} converges to the solution of \eqref{target1}. The convergence seems generally fast as shows the example of Table~\ref{tableConv}.
\begin{table}
$$
\begin{array}{|c|c|c|c|c|c|c|c|c|c|}
\hline
n & 0 & 1 & 2 & 3 & 4 & 5 & 6 & 7 & 8\\
\hline
\hline
\multirow{2}{*}{$M_n^{(1)}$} & 5.95 & 5.95 & 5.95 & 5.95 & 5.95 & 5.95 & 5.95 & 5.95 & 5.95 \\
& 2.62 & 4.70 & 4.09 & 4.23 & 4.19  & 4.20 & 4.20 & 4.20 & 4.20 \\
\hline
\textbf{$p_i^{(1)}$} &\textbf{ 3.64 }& \textbf{5.25} & \textbf{4.85} & \textbf{4.94} & \textbf{4.91} & \textbf{4.92} & \textbf{4.92} & \textbf{4.92} & \textbf{4.92}\\
\hline
\hline
\multirow{2}{*}{$M_n^{(2)}$} & 6.03 & 3.36 & 3.86 & 3.74 & 3.77 & 3.76 & 3.76 & 3.76 & 3.76 \\
& 7.11 & 7.11 & 7.11 & 7.11 & 7.11 & 7.11 & 7.11 & 7.11 & 7.11 \\
\hline
\textbf{$p_i^{(2)}$} & \textbf{6.53} & \textbf{4.57} & \textbf{5.01} & \textbf{4.90} & \textbf{4.93} & \textbf{4.92} & \textbf{4.92} & \textbf{4.92} & \textbf{4.92}\\
\hline
\end{array}
$$
\label{tableConv}
\caption{Realization of \eqref{algo} with $K=1$, $L=2$ and $f_1=f_2=1$, given $M_0^{(1)}$ and $M_0^{(2)}$.\label{tableConv}}
\end{table}
This yields
\begin{equation} \label{alphasol}
\alpha^{(l)}= \alpha\frac{f^{(l)}}{p_{0}^{(l)}}\frac{f^{(l)}}{p_{1}^{(l)}}\frac{f^{(l)}}{p_{l}^{(l)}}\cdots\qquad (1\leqslant l \leqslant L).
\end{equation}
where any $\alpha>0$ suits to obtain \eqref{target1}.
To simultaneously obtain \eqref{target2}, we know by \ref{relfreq} that the solution of \eqref{target1} and \eqref{target2} together exists and is unique, thus, as $n\to\infty$ the value of
$\sum_{j\in\Z}\sum_{i\geqslant0} a_{i,n}^{(j)}(x)$ cannot depend on $x$.
Although, because of numerical reasons, it may slightly vary with $x$.
Thus we suggest to keep the dependence in $x$ and replace $\alpha$ in \eqref{alphasol} by
$$ \alpha(x) = \lim_{n\to\infty} \frac{\displaystyle 1}{ \sum_{j\in\Z}\sum_{i\geqslant0} a_{i,n}^{(j)}(x)} $$
in \eqref{alphasol}.
This completes the procedure to estimate the $a_i^{(j)}$ form observations.

\begin{figure}[here!]
\hspace{-1.6cm}\includegraphics[scale=0.35]{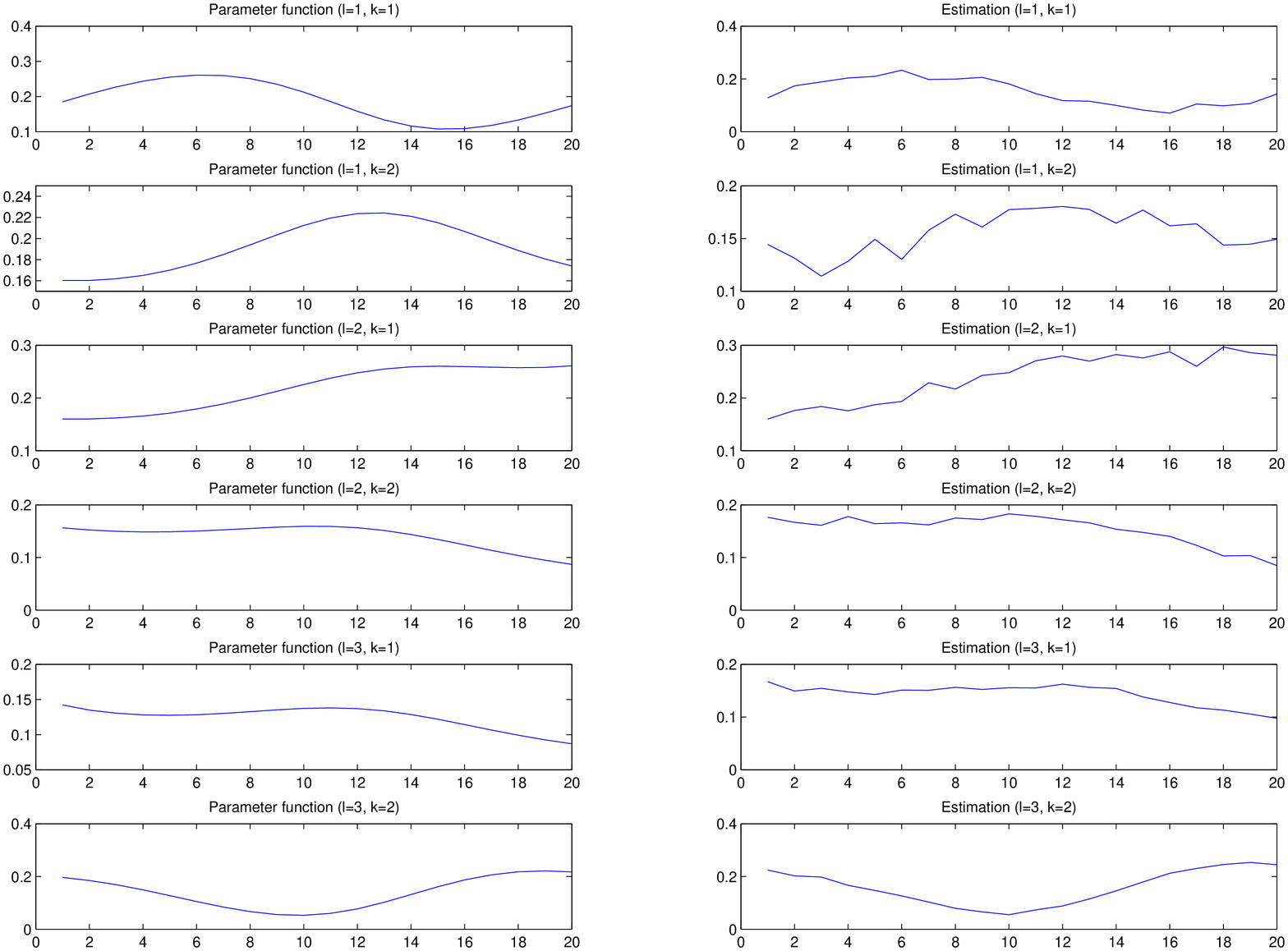}
\begin{pspicture}(-7,0)(7,0)
\psset{xunit=1cm,yunit=1cm,runit=1cm}
\uput[u](-1.5,.2){\footnotesize{$T=500$}}
\end{pspicture}
\caption{Recovery of the parameter functions with $D=20$, $K=2$, $L=3$, $\sigma=1$ and PAM.}
\label{EstFig}
\end{figure}

Figure~\ref{EstFig} shows an output of the full algorithm using PAM, with measurement errors but given the true values of $K$ and $L$.

\subsection{Distance between two sets of parameter functions} \label{hDist}

To measure the quality of the estimation, it is necessary to quantify the dissimilarity between the the original parameter functions $a_i^{(j)}$ ($0\leqslant i < K$, $1\leqslant j \leqslant L$) and their estimations $\hat{a}_i^{(j)}$ ($0\leqslant i < \hat{K}$, $1\leqslant j \leqslant \hat{L}$). If $\hat{K}\neq K$ the estimation can certainly be qualified as bad, so that only the case $\hat{K}=K$ requires a discussion.

The order in which the different patterns are retrieved can change and their total numbers can differ. Consequently the Hausdorff distance  between the $L$ graphs of $(a_i^{(j)}(x))_{i}$  ($1\leqslant j \leqslant L$) and the $\hat{L}$ graphs of $(\hat{a}_i^{(j)}(x))_{i}$ ($1\leqslant j \leqslant \hat{L}$), that are compact in $[0,1]^q\times\R^K$ (or in $\{1,\ldots,D\}\times\R^K\cong\R^{DK}$ for the discrete version), perfectly suits.
Recall that the Hausdorff distance between nonempty compacts $A$ and $B$ is
$$ d_{\mathcal{H}}(A,B)=\max\{\sup_{a\in A}d(a,B),\sup_{b\in b}d(b,A)\}, $$
here considered with the Euclidean distance.
We have thus to compute $d_{\mathcal{H}}(A,B)$ with
$$ A=\bigcup_{j=1}^L \operatorname{graph} [(a_i^{(j)}(x))_{0\leqslant i<K}] \qquad\textrm{and}\qquad B=\bigcup_{j=1}^{\hat{L}} \operatorname{graph} [(\hat{a}_i^{(j)}(x))_{0\leqslant i<K}]$$
to reach the stated goal.

\begin{figure}[here!]
\begin{center}
\includegraphics[scale=0.5]{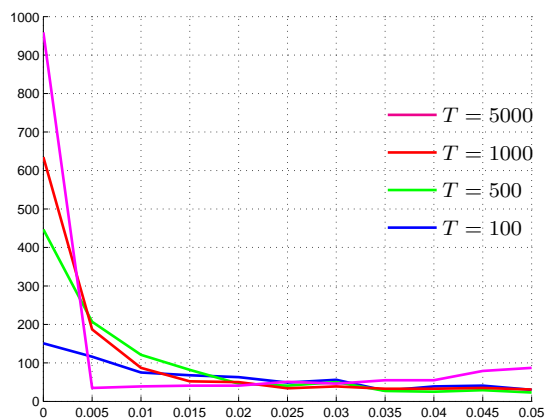}
\begin{pspicture}(-7,0)(7,0)
\psset{xunit=1cm,yunit=1cm,runit=1cm}
\uput[r](1.1,3.35){\footnotesize{$T=100$}}
\uput[r](1.1,3.85){\footnotesize{$T=500$}}
\uput[r](1.1,4.35){\footnotesize{$T=1000$}}
\uput[r](1.1,4.85){\footnotesize{$T=5000$}}
\psline[linewidth=1pt,plotstyle=curve,linecolor=blue]{-}(.55,3.35)(1.2,3.35)
\psline[linewidth=1pt,plotstyle=curve,linecolor=green]{-}(.55,3.85)(1.2,3.85)
\psline[linewidth=1pt,plotstyle=curve,linecolor=red]{-}(.55,4.35)(1.2,4.35)
\psline[linewidth=1pt,plotstyle=curve,linecolor=magenta]{-}(.55,4.85)(1.2,4.85)

\end{pspicture}
\caption{Smoothed histograms of the Hausdorff distance between the original parameter functions $(a_i^{(j)}(x))_{0\leqslant i<K}$ ($1\leqslant j \leqslant L$) and their estimations $(\hat{a}_i^{(j)}(x))_{0\leqslant i<K}$ ($1\leqslant j \leqslant \hat{L}$).}
\label{AsBeh}
\end{center}
\end{figure}

To illustrate the procedure developed in \ref{estA}, Figure~\ref{AsBeh} shows smoothed histograms of the distances from \ref{hDist} for the estimation of the parameter functions with sample sizes $T=100$, $500$, $1000$, $5000$ given the true values of $K$ and $L$. The test was performed with $N=5000$ trials for each sample size: $C$ running in $\{2,4,6,8,10\}$, $D$ in $\{1,5,10,15,20\}$, $K$ from 2 to 5, $L$ from 1 to 5 with 10 repetitions of each. The other parameters were $Q=\min(\operatorname{ceil}(\underline{p}T),100)$ with $\underline{p}$ from \ref{relfreq} and $\sigma=1$.

\section{Conclusion}

After the linear processes in function spaces with heavy tailed innovations \cite{MS10}, CM3 processes are other examples of jointly regularly varying time series in function spaces. Under \eqref{finiteSum} they enjoy the finite-cluster condition and under \eqref{fastConv} the strong mixing property.

Further studies could be to determine whether or not the approximation theorem of Deheuvels for M3 \cite{DEH83} and of Smith and Weissman for M4 \cite{SW96} also hold for CM3, that is to say if max-stable processes in function spaces, excluding the ones containing a deterministic component, can be arbitrary closely approached by a CM3. From these papers also arises the question of a generalization of the multivariate extremal index. Since such an object becomes hard to figure out in function spaces, it has maybe to be replaced by the spectral process.

About the estimation, in the empirical study on CM3, we saw that the mode correctly estimates more frequently the length of the tail dependence than the median and the mean. For finding the number of patterns, this study revealed the importance of simulating the behavior of the chosen method before using it on real data.

\section*{Acknowledgements}
The author thanks Johan Segers for helpful discussions throughout the writing of this paper and especially for his written
communications about the finite-cluster condition and the strong mixing condition.

\vspace{1cm}

\noindent The Matlab code written for this paper is available on the Matlab Central File Exchange at \href{http://www.mathworks.com/matlabcentral/fileexchange/26539}{http://www.mathworks.com/matlabcentral/fileexchange/26539}.

\newpage

\appendix

\section{Proofs of Section \ref{SectionCM3}} \label{proofs}

\noindent\textbf{Proof of Proposition~\ref{CM3processLem}.}
\textsc{Claim 1} - \textit{$X_t$ is a random element in $\mathcal{C}([0,1]^q,\R_+)$.}
Given $t\in\Z$, write $X\overset{d}{=}X_t$ and $X^{(m)} = \sup_{i,|j| \le m} a_i^{(j)} Z_{t-i}^{(j)}.$
For $z > 0$, since the $Z_i^{(j)}$ are iid unit-Fr\'{e}chet,
\begin{align*}
    P( \| X \|_\infty \le z )
    &= P( \forall i\geqslant0,\forall j \in\Z : Z_{t-i}^{(j)} \le z / \| a_i^{(j)} \|_\infty ) \\
    &= \biindice{\displaystyle\prod}{j \in\Z}{i\geqslant0} \exp \bigl( - \| a_i^{(j)} \|_\infty / z \bigr)
    = \exp \bigl( - {\textstyle\sum_{j \in\Z}\sum_{i\geqslant0}} \| a_i^{(j)} \|_\infty / z \bigr).
\end{align*}
If $\sum_{j \in\Z}\sum_{i\geqslant0} \| a_i^{(j)} \|_\infty = \infty$, then this probability is zero for all $z > 0$,
hence $\| X \|_\infty = \infty$ with probability one.

If $\sum_{j \in\Z}\sum_{i\geqslant0} \| a_i^{(j)} \|_\infty < \infty$, then
$
    X^{(m)} \le X \le X^{(m)} + \sup_{i,|j| > m} \| a_i^{(j)} \|_\infty Z_{t-i}^{(j)}
$
and thus
$
    \| X - X^{(m)} \|_\infty \le \sup_{i,|j| > m} \| a_i^{(j)} \|_\infty Z_{t-i}^{(j)}.
$
Writing $$Y_m = \sup_{i,|j| > m} \| a_i^{(j)} \|_\infty Z_i,$$ we have $Y_m \ge 0$ and $Y_m$ is decreasing in $m$.
By a similar computation as above, we find that for $y > 0$,
$$
    P(Y_m \le y) = \exp \bigl( - {\textstyle\sum_{i,|j| > m}} \| a_i^{(j)} \|_\infty / y \bigr).
$$
As a consequence, $Y_m \to 0$ in probability and, by monotonicity, almost surely.
Since the uniform limit of a sequence of continuous functions is continuous, by monotonicity, $X$ is continuous with probability one.

Then the fact that, for every $t\in \Z$, the map $\omega\mapsto X_t(\omega,\cdot)$ with values in $\mathcal{C}([0,1]^q,\R_+)$ is measurable follows easily.

\textsc{Claim 2} - \textit{$(X_t)_{t\in Z}$ is a stationary time series.}
In extension, this property is that, for every $n\geqslant0$, every $h\geqslant0$ and every measurable set $A\subset\mathcal{C}([0,1]^q,\R^h)$,
$$ P((X_0,\ldots,X_n)\in A)=P((X_{h},\ldots,X_{n+h})\in A). $$
The argument is based on two facts:\\
1) It suffices to pick up $x\in[0,1]^q$, $z_0,\ldots,z_n\in\R$ and verify that $$ P(X_0(x)\leqslant z_0,\ldots,X_n(x)\leqslant z_n)
=P(X_h(x)\leqslant z_0,\ldots,X_{n+h}(x)\leqslant z_n), $$
the case with $k$ points $x_1,\ldots x_k\in[0,1]^q$ being similar.

\noindent 2) The right hand-side of the expression in 1) admits the expansion
$$\begin{array}{l}
\displaystyle P(X_h(x)\leqslant z_0,\ldots,X_{n+h}(x)\leqslant z_n)\\
\rule{0cm}{.4cm} \displaystyle = P(\sup_{j\in\Z}\sup_{i\geqslant0}a_i^{(j)}(x)Z_{h-i}^{(j)}\leqslant z_0,\ldots,\sup_{j\in\Z}\sup_{i\geqslant0}a_i^{(j)}(x)Z_{n+h-i}^{(j)}\leqslant z_n)\\
\rule{0cm}{.4cm} \displaystyle = P(\forall j_0\in\Z,\forall i_0\geqslant0:a_{i_0}^{(j_0)}(x)Z_{h-i_0}^{(j_0)}\leqslant z_0,\ldots,\\
\rule{0cm}{.5cm} \displaystyle  \qquad\qquad\qquad\qquad\qquad\qquad\qquad  \forall j_n\in\Z, \forall i_n\geqslant0:a_{i_n}^{(j_n)}(x)Z_{n+h-i_n}^{(j_n)}\leqslant z_n)\\
\rule{0cm}{.4cm} \displaystyle = P(\forall j\in\Z: Z_{h+n}^{(j)}\leqslant\frac{z_n}{a_0^{(j)}(x)},\\
\rule{0cm}{.4cm} \displaystyle  \qquad\qquad\qquad \ Z_{h+n-1}^{(j)}\leqslant\min\left(\frac{z_n}{a_1^{(j)}(x)},\frac{z_{n-1}}{a_0^{(j)}(x)}\right),\\
\rule{0cm}{.4cm} \displaystyle  \qquad\qquad\qquad \ Z_{h+n-2}^{(j)}\leqslant\min\left(\frac{z_n}{a_2^{(j)}(x)},\frac{z_{n-1}}{a_1^{(j)}(x)},\frac{z_{n-2}}{a_0^{(j)}(x)}\right),\\
\rule{0cm}{.4cm} \displaystyle  \qquad\qquad\qquad \ \ldots,\\
\rule{0cm}{.4cm} \displaystyle  \qquad\qquad\qquad \ Z_{h}^{(j)}\leqslant\min\left(\frac{z_n}{a_n^{(j)}(x)},\frac{z_{n-1}}{a_{n-1}^{(j)}(x)},\ldots,\frac{z_{0}}{a_0^{(j)}(x)}\right),\\
\rule{0cm}{.4cm} \displaystyle  \qquad\qquad\qquad \ Z_{h-1}^{(j)}\leqslant\min\left(\frac{z_n}{a_{n+1}^{(j)}(x)},\frac{z_{n-1}}{a_{n}^{(j)}(x)},\ldots,\frac{z_{0}}{a_1^{(j)}(x)}\right),\\
\rule{0cm}{.4cm} \displaystyle  \qquad\qquad\qquad \ \ldots)\\
\end{array}$$
which allows us to dispose of the variable $h$ by independence.\\

\noindent\textbf{Proof of Proposition~\ref{CM3process}.}
Given $s,t\geqslant0$ consider the Dirac masses in $\mathcal{C}([0,1]^q,\R)$
$$\begin{array}{rcl}
\displaystyle\delta(-s,t;i,j)&\overset{d}{=}&\displaystyle \left(\frac{a_{-s+i}^{(j)}}{\| a_{i}^{(j)} \|_\infty} ,\ldots, \frac{a_{t+i}^{(j)}}{\| a_i^{(j)} \|_\infty}\right)\\
\rule{0cm}{.7cm}&\overset{d}{=}&\displaystyle \left[ \left.\left(\frac{a_{-s+i}^{(j)}Z}{\| a_{i}^{(j)}Z \|_\infty} ,\ldots, \frac{a_{t+i}^{(j)}Z}{\| a_i^{(j)}Z \|_\infty}\right) \ \right|\ \| a_i^{(j)} \|_\infty Z> x \right]\\
\end{array}$$
with $i\geqslant0$ and $j\in\Z$.

A slight variation in the proof of \cite[Lemma B.2]{MS10} using the identity
\begin{multline*}
     1\{ \max( \norm{X_1}, \norm{X_2}) > x \}
    = \\ 1\{\norm{X_1} > x\} + 1\{\norm{X_2} > x\} - 1\{\norm{X_1} > x\} 1\{\norm{X_2} > x\}
\end{multline*}
and considering $T_i^{(j)}(z)=a_i^{(j)}z$ leads to the new conclusion that
\begin{equation} \label{supmax}
    \E \biggl| 1{\{\sup_{j\in\Z}\sup_{i\geqslant0} \norm{a_i^{(j)}}Z_{-i}^{(j)} > x \}}
    - \biindice{\displaystyle\sum}{j \in\Z}{i\geqslant0}  1{\{\norm{a_i^{(j)}}Z > x) \}} \biggr|= o \bigl( P(Z>x) \bigr),\, x \to \infty.
\end{equation}
Consequently, for $f\in\mathcal{C}_{b}(\mathcal{C}([0,1]^q,\R^{s+t+1}))$, by virtue of \eqref{supmax}
$$\begin{array}{l}
\displaystyle E\left[\left.  f\left(\frac{X_{-s}}{x},\ldots,\frac{X_t}{x}\right)  \ \right|\ \|X_0\|_\infty>x  \right]\\
\rule{0cm}{.9cm}\displaystyle =\frac{E\left[ f\left(\frac{X_{-s}}{x},\ldots,\frac{X_{t}}{x}\right)  1\{\|X_0\|_\infty>x\}  \right]}{P(Z>x)} \frac{P(Z>x)}{P(\|X_0\|_\infty>x)}\\
\rule{0cm}{.9cm}\displaystyle =\frac{E\left[ f\left(\frac{X_{-s}}{x},\ldots,\frac{X_{t}}{x}\right)  1\{\displaystyle\sup_{j\in\Z}\sup_{i\geqslant0} \|a_i^{(j)}\|_\infty Z_{-i}^{(j)}>x\}  \right]}{P(Z>x)} \frac{P(Z>x)}{P(\|X_0\|_\infty>x)}\\
\rule{0cm}{.9cm}\displaystyle = \biindice{\displaystyle\sum}{j \in\Z}{i\geqslant0} \frac{E\left[ f\left(\frac{X_{-s}}{x},\ldots,\frac{X_{t}}{x}\right)  1\{\|a_i^{(j)}\|_\infty Z_{-i}^{(j)}>x\}  \right]}{P(\|a_i^{(j)}\|_\infty Z_{-i}^{(j)}>x)} \\
\displaystyle \hspace{1.8cm} \cdot \frac{P(\|a_i^{(j)}\|_\infty Z_{-i}>x)}{P( Z>x)} \frac{P(Z>x)}{P(\|X_0\|_\infty>x)} +o(1).\\
\end{array}$$
Thanks to the dominated convergence and the continuity of $f$, the last expression is equal to
$$\begin{array}{l}
\displaystyle \biindice{\displaystyle\sum}{j \in\Z}{i\geqslant0} E\left[\left. f\left(\frac{\displaystyle\sup_{l\in \Z}\sup_{k \geqslant 0} (a_k^{(l)} Z_{-s-k}^{(l)})}{x},\ldots,\frac{\displaystyle\sup_{l\in \Z}\sup_{k \geqslant 0} (a_k^{(l)} Z_{t-k}^{(l)})}{x}\right) \right| \|a_i^{(j)}\|_\infty Z_{-i}^{(j)}\!>\!x  \right] \\
\displaystyle \hspace{\stretch{1}} \cdot \frac{P(\|a_i^{(j)}\|_\infty Z>x)}{P( Z>x)} \frac{P(Z>x)}{P(\|X_0\|_\infty>x)} +o(1)\\
\displaystyle =\biindice{\displaystyle\sum}{j \in\Z}{i\geqslant0} E\left[\left. f\left(\frac{a_{-s+i}^{(j)} Z}{x} ,\ldots, \frac{a_{t+i}^{(j)} Z}{x}\right) \ \right|  \|a_i^{(j)}\|_\infty Z>x  \right]\\
\displaystyle \hspace{\stretch{1}} \cdot  \frac{P(\|a_i^{(j)}\|_\infty Z>x)}{P( Z>x)} \frac{P(Z>x)}{P(\|X_0\|_\infty>x)} +o(1)\\
\displaystyle =\biindice{\displaystyle\sum}{j \in\Z}{i\geqslant0} E\left[\left. f\left(\frac{Z}{x}(a_{-s+i}^{(j)},\ldots, a_{t+i}^{(j)}\right) \ \right|  \|a_i^{(j)}\|_\infty Z>x  \right]\\
\displaystyle \hspace{\stretch{1}} \cdot
\frac{P(\|a_i^{(j)}\|_\infty Z>x)}{P( Z>x)} \frac{P(Z>x)}{P(\|X_0\|_\infty>x)} +o(1).\\
\end{array}$$
Using the continuous mapping theorem and the regular variation of $Z$, the last term converges to
$$\begin{array}{l}
\displaystyle \biindice{\displaystyle\sum}{j \in\Z}{i\geqslant0} E\left[f(Y\delta(-s,t;i,j))\right] \|a_i^{(j)}\|_\infty \frac{1}{\sum_{l\in\Z}\sum_{k\geqslant0} \|a_k^{(l)}\|_\infty}\\
\end{array}$$
where $Y\sim\Pareto(1)$. The factors
$$ \pi_{ij}:=  \frac{\|a_i^{(j)}\|_\infty}{\sum_{l\in\Z}\sum_{k\geqslant 0} \|a_k^{(l)}\|_\infty} $$
form a probability distribution on $\Z_+\times\Z$, let us say of a random vector $S$.

According to \cite[Theorem 3.1 (iv)]{MS10}, the spectral process of the CM3 process $(X_t)_{t\in Z}$ is
$$ (\Theta_{-s},\ldots,\Theta_{t}) \overset{d}{=} \delta(-s,t;S) $$
which has the announced distribution.\\

\noindent\textbf{Proof of Proposition~\ref{conditionCThm}.}
It suffices to check that
$$\lim_{m\to\infty} \limsup_{n\to\infty} P(\exists t \in\{2m,\ldots,r_n\} : \|X_t\|>n\ |\ \|X_0\|>n )=0.$$
For the convenience of the proof, set $p=t-i$ so that
$$ X_t= \sup_{j\in\Z} \sup_{p\leqslant t} a_{t-p}^{(j)} Z_{p}^{(j)}. $$
Let $b_i^{(j)}:=\|a_i^{(j)}\|_\infty$ to have
$$ \|X_t\|= \sup_{j\in\Z} \sup_{p\leqslant t} b_{t-p}^{(j)} Z_{p}^{(j)}. $$
Given $m>0$, write $A_n:=\{\|X_0\|>n\}$ so that
$$
\begin{array}{l}
A_{n}=\{ \exists j\in\Z,\exists p\leqslant 0 : b_{-p}^{(j)}Z_{p}^{(j)}>n  \}.\\
\end{array}
$$
We have, as $n\to\infty$,
$$\begin{array}{l}
nP(A_{n})\\
\rule{0cm}{.7cm}=n[1-P(\forall j\in\Z,\forall p\leqslant 0 : b_{-p}^{(j)}Z_{p}^{(j)} \leqslant n)]\\
\rule{0cm}{.7cm} =n[1-\exp(-\frac{1}{n} \sum_{j\in\Z} \sum_{p\leqslant 0} b_{-p}^{(j)})]\\
\rule{0cm}{.7cm}\to\sum_{j\in\Z} \sum_{p\leqslant 0} b_{-p}^{(j)} = \sum_{j\in\Z} \sum_{p\geqslant 0} b_{p}^{(j)}<+\infty. \\
\end{array}$$
For $m>0$, decompose $\{\exists t \in\{2m,\ldots,r_n\} : \|X_t\|>n\}=C_{m,n}\dot{\cup} D_{m,n}$ where
$$
\begin{array}{l}
C_{m,n}=\{ \exists t \in\{2m,\ldots,r_n\},\exists j\in\Z,\exists p> m : b_{t-p}^{(j)}Z_{p}^{(j)}>n  \},\\
D_{m,n}=\{ \exists t \in\{2m,\ldots,r_n\},\exists j\in\Z,\exists p\leqslant m : b_{t-p}^{(j)}Z_{p}^{(j)}> n  \}.\\
\end{array}
$$
We have
$$\begin{array}{rcl}
P(C_{m,n})&=&\displaystyle 1-P(\forall t \in\{2m,\ldots,r_n\},\forall j\in\Z,\forall p> m : b_{t-p}^{(j)}Z_{p}^{(j)} \leqslant n)\\
\rule{0cm}{.4cm}&=&\displaystyle 1-P(\forall j\in\Z,\forall p> m : [\max_{2m\leqslant t \leqslant r_n} b_{t-p}^{(j)}]Z_{p}^{(j)} \leqslant n)\\
\rule{0cm}{.4cm}&=&\displaystyle 1-\exp(-\frac{1}{n} \sum_{j\in\Z} \sum_{p> m} \max_{2m\leqslant t \leqslant r_n} b_{t-p}^{(j)}) \\
\rule{0cm}{.4cm}&\leqslant&\displaystyle \frac{1}{n} \sum_{j\in\Z} \sum_{p> m} \max_{2m\leqslant t \leqslant r_n} b_{t-p}^{(j)} \\
\rule{0cm}{.4cm}&\leqslant&\displaystyle \frac{m_n}{n} \sum_{j\in\Z} \sum_{k\geqslant0} b_{k}^{(j)} \to 0 \\
\end{array}$$
as $n\to\infty$.

Setting $p=2m-q$ write
$$ A_{n}=\{ \exists j\in\Z,\exists q\geqslant 2m : b_{q-2m}^{(j)}Z_{q}^{(j)}>n  \} $$
and
$$ D_{m,n}=\{ \exists j\in\Z,\exists q\geqslant m : \max_{q\leqslant t \leqslant q-2m + r_n} b_{t}^{(j)}Z_{q}^{(j)}> n  \} $$
hence
$$\begin{array}{rl}
nP(A_{n}\cap D_{m,n})&\displaystyle\leqslant n\sum_{j\in\Z} \sum_{q\geqslant 2m} P(\{ b_{q-2m}^{(j)}Z_{q}^{(j)}>n \}\cap D_{m,n})\\
\rule{0cm}{.7cm}&\displaystyle\leqslant n\sum_{j\in\Z} \sum_{q\geqslant 2m} P(\min[ b_{q-2m}^{(j)},\max_{t\geqslant q } b_{t}^{(j)}]Z_{q}^{(j)} >n) \\
\rule{0cm}{.7cm}&\displaystyle\hspace{.3cm}+ n\sum_{j\in\Z} \sum_{q\geqslant 2m} P( b_{q-2m}^{(j)}Z_{q}^{(j)}>n ) P(D_{m,n}). \\
\end{array}$$
The first term of the sum is bounded above by
$$\begin{array}{l}
\displaystyle\sum_{j\in\Z} \sum_{q\geqslant 0}  \min[ b_{q}^{(j)},\max_{t\geqslant2 m } b_{t}^{(j)}]\\
\end{array}$$
and the second term tends to $0$ as $n\to\infty$ since
$$
 n \sum_{j\in\Z} \sum_{q\geqslant 2m}  P( b_{q-2m}^{(j)}Z_{q}^{(j)}>n )\leqslant\sum_{j\in\Z} \sum_{q\geqslant 0} b_{q}^{(j)}  <\infty
$$
and
$$
P(D_{m,n})= \frac{1}{n}\sum_{j\in\Z} \sum_{q\geqslant m} \max_{q\leqslant t \leqslant q-2m + r_n} b_{t}^{(j)}
\leqslant\frac{m_n}{n}\sum_{j\in\Z} \sum_{q\geqslant0} b_{q}^{(j)} \to 0.
$$
Consequently
$$ \limsup_{n\to\infty}nP(A_{n}\cap D_{m,n})\leqslant \sum_{j\in\Z} \sum_{p\in\Z} \min[ b_{p}^{(j)},\max_{t\geqslant 2m } b_{t}^{(j)}].$$

Since $A_{n}$ and $C_{m,n}$ are independent,
$$\begin{array}{l}
\displaystyle P(\exists t \in\{2m,\ldots,r_n\} : \|X_t\|>n\ |\ \|X_0\|>n )\\
\rule{0cm}{.7cm} \displaystyle = P(C_{m,n} \cup D_{m,n}\ |\ A_{n} )\\
\rule{0cm}{.7cm} \displaystyle = \frac{P([C_{m,n} \cup D_{m,n}] \cap A_{m,n} )}{P(A_{n}) }\\
\rule{0cm}{.7cm} \displaystyle \leqslant \frac{P([C_{m,n} \cap A_{n}] \cup [ D_{m,n}\cap A_{n} ]  \}) }{P(A_{n} )}\\
\rule{0cm}{.7cm} \displaystyle \leqslant P(C_{m,n} ) + \frac{P( D_{m,n} \cap A_{n} )}{P(A_{n}) }\\
\end{array}$$
so that
$$\begin{array}{l}
\displaystyle \limsup_{n\to\infty} P(\exists t \in\{2m,\ldots,r_n\} : \|X_t\|>n\ |\ \|X_0\|>n )\\
\rule{0cm}{.7cm} \displaystyle \leqslant \frac{1}{\displaystyle\sum_{j\in\Z}\sum_{p\geqslant 0} b_{p}^{(j)}}\left[\sum_{j\in\Z} \sum_{p\in\Z} \min[ b_{p}^{(j)},\max_{t\geqslant 2m } b_{t}^{(j)}] \right]\\
\end{array}$$
which tends to $0$ if $m\to\infty$.\\

\noindent\textbf{Proof of Proposition~\ref{conditionSMThm}.}
Again set $p=t-i$ to have
$$ X_t= \sup_{j\in\Z} \sup_{p\leqslant t} a_{t-p}^{(j)} Z_{p}^{(j)}. $$
For large $t\geqslant0$, we will prove that we can approach $X_t$ by
$$ X_t^{(+)}= \sup_{j\in\Z} \sup_{1\leqslant p\leqslant t} a_{t-p}^{(j)} Z_{p}^{(j)} $$
in the sense of $$ \lim_{m\to+\infty}P(\forall t\geqslant m : X_{t}^{(+)}=X_t)=1 .$$
Then the conclusion will follow from the fact that the processes $(X_t^{(+)})_{t\geqslant1}$ and $(X_t)_{t\leqslant-1}$ are independent.

As the Borel $\sigma$-field on $\mathcal{C}([0,1]^q)$ is generated by the finite dimensional sets, it is sufficient to check that
the limit is 1 on every subset $\{x_1,\ldots,x_k\} \subset [0,1]^q$. Passing to the complementary, compute
$$\begin{array}{rcl}
p_m&:=&\displaystyle P(\exists t \geqslant m,\exists 1\leqslant i \leqslant k \, :\, X_t^{(+)}(x_i) < X_t(x_i))\\
\rule{0cm}{.4cm}&=&\displaystyle P(\exists t \geqslant m,\exists 1\leqslant i \leqslant k \, :\, \sup_{j\in\Z} \sup_{1\leqslant p\leqslant t} a_{t-p}^{(j)}(x_i) Z_{p}^{(j)}< \\
&&\hspace{5cm}  \sup_{j\in\Z} \sup_{p\leqslant t} a_{t-p}^{(j)}(x_i) Z_{p}^{(j)})\\
\rule{0cm}{.4cm}&=&\displaystyle P(\exists t \geqslant m,\exists 1\leqslant i \leqslant k \, :\, \sup_{j\in\Z} \sup_{1\leqslant p\leqslant t}  a_{t-p}^{(j)}(x_i) Z_{p}^{(j)}< \\
&&\hspace{5cm} \sup_{j\in\Z} \sup_{p\leqslant 0} a_{t-p}^{(j)}(x_i) Z_{p}^{(j)}).\\
\end{array}$$
If $a,b>0$ and $Z_1,Z_2$ are standard unit-Fr\'echet, as in \eqref{profprob},
$$
P(aZ_1<bZ_2)=\frac{1}{1+\frac{a}{b}}<\frac{b}{a}
$$
so that
$$ p_m \leqslant \sum_{t \geqslant m}\sum_{i=1}^k \frac{\sum_{j\in\Z} \sum_{p\leqslant 0} a_{t-p}^{(j)}(x_i)}{\sum_{j\in\Z} \sum_{1\leqslant p\leqslant t}  a_{t-p}^{(j)}(x_i)}=:\sum_{t \geqslant m}\sum_{i=1}^k \frac{N_t(x_i)}{D_t(x_i)}. $$
According to \eqref{finiteSum}, $D_t(x)$ uniformly converges to a continuous function $D(x)$ as $t\to\infty$.
Without loss of generality, we can assume that $D(t)\equiv 1$.
Write $D_t(x)=1-\varepsilon_t(x)$ with $\|\varepsilon_t\|_\infty\to0$ as $t\to\infty$. Consequently,
$$ \frac{1}{D_t(x)}=1+\frac{\varepsilon_t(x)}{D_t(x)} $$
hence
$$ p_m \leqslant \sum_{t \geqslant m}\sum_{i=1}^k N_t(x_i)  + \sum_{t \geqslant m}\sum_{i=1}^k N_t(x_i) \frac{\varepsilon_t(x_i)}{D_t(x_i)}.$$
Next, use the fact that for every $\eta>0$, there exists $m_\eta$ such that for every $t\geqslant m_\eta$ and every $1\leqslant i \leqslant k$
$$ 0 \leqslant \frac{\varepsilon_t(x_i)}{D_t(x_i)} \leqslant \eta $$
to see that $$ \frac{\displaystyle\sum_{t \geqslant m}\sum_{i=1}^k N_t(x_i) \frac{\varepsilon_t(x_i)}{D_t(x_i)}}{\displaystyle\sum_{t \geqslant m}\sum_{i=1}^k N_t(x_i) } = o(1) $$
and so obtain
$$ p_m \leqslant (1+o(1)) \sum_{t \geqslant m}\sum_{i=1}^k \sum_{j\in\Z} \sum_{p\leqslant 0} a_{t-p}^{(j)}(x_i)  .$$
Since $$ \sum_{t \geqslant m} \sum_{i=1}^k \sum_{j\in\Z} \sum_{p\leqslant 0} a_{t-p}^{(j)}(x_i) = \sum_{k \geqslant m}\sum_{i=1}^k \sum_{j\in\Z} (k-m+1) a_{k}^{(j)}(x_i),$$
thanks to \eqref{fastConv} we have that $p_m$ vanishes as $m\to\infty$, which proves the result.\\

\noindent\textbf{Proof of Proposition~\ref{extIndex}.}
Since $$ \|X_t\|_\infty = \|\sup_{j\in\Z} \sup_{i\geqslant0} a_i^{(j)} Z_{t-i}^{(j)}\|_\infty
=  \sup_{j\in\Z} \sup_{i\geqslant0} \|a_i^{(j)}\|_\infty Z_{t-i}^{(j)}$$
it suffices to transcribe the formula for a M3 from \cite{SW96}.

\end{document}